\documentclass[iicol, sn-mathphys]{sn-jnl}
\usepackage[utf8]{inputenc}
\usepackage{enumitem}
\usepackage{soul, ulem}
\usepackage{yfonts}
\usepackage[letterspace]{microtype}
\usepackage[export]{adjustbox}
\usepackage{anyfontsize}
\usepackage{makecell}
\usepackage[caption=false]{subfig}
\usepackage{multirow}%
\usepackage{amsfonts}%
\usepackage[title]{appendix}%
\usepackage{xcolor}%
\usepackage{textcomp}%
\usepackage{manyfoot}%
\usepackage{listings}%
\setcitestyle{numbers}
\setcitestyle{square}

\RequirePackage[OT1]{fontenc}

\RequirePackage{amsmath}
\RequirePackage{hyperref}

\usepackage[capitalise]{cleveref}
\usepackage{algorithm}
\usepackage{algorithmic}
\usepackage{amssymb,amsbsy}
\usepackage{amscd}
\usepackage{latexsym}
\usepackage{vmargin,fancyhdr}
\usepackage{graphicx}
\usepackage{dsfont}
\usepackage{booktabs, siunitx}
\usepackage{color}
\usepackage{mathtools}

\newenvironment{simpletabular}{\setlength{\extrarowheight}{1pt}\setlength{\tabcolsep}{1.5pt}
  \begin{tabular}}{\end{tabular}}

\newcommand{\rev}[1]{{\color{black} #1}}

\DeclarePairedDelimiter\ceil{\lceil}{\rceil}

\newcommand{\Z}{\mathbb{Z}}

\newcommand{\Kg}{K_g}
\newcommand{\scs}{\scriptstyle}
\newcommand{\closer}[1]{\textls[-100] {#1}}
\newcommand{\close}[1]{\textls[-50] {#1}}
\newcommand{\DeltaJROI}{{$\Delta_{J}(\text{\footnotesize ROI}) \,$}}
\newcommand{\smax}{S_{max}}

\newcommand{\tableLDDMMmultiscale}{\makecell{\textbf{LDDMM} \\ $+$ \textbf{Multiscale}}}
\newcommand{\tableLDDMMmultiscalee}{\makecell{\closer{LDDMM}$+$\\ \close{Multiscale}}}

\def\locpix{x}
\def\gk{c_{k}}

\newcommand{\To}{\longrightarrow}
\newcommand{\flowphi}{\Phi}

\begin{document}
\title[Article Title]{Wavelet-Based Multiscale Flow For Realistic Image Deformation in the Large Diffeomorphic Deformation Model Framework}

\author*[1]{\fnm{Fleur} \sur{Gaudfernau}}\email{fleur.gaudfernau@etu.u-paris.fr}

\author[2]{\fnm{El\'eonore} \sur{Blondiaux}}\email{eleonore.blondiaux@aphp.fr}

\author[1]{\fnm{St\'ephanie} \sur{Allassonni\`ere}}\email{stephanie.allassonniere@u-paris.fr}

\author[3]{\fnm{Erwan} \sur{Le Pennec}}\email{erwan.le-pennec@polytechnique.edu}

\affil[1]{\orgdiv{CRC, Université de Paris, INRIA EPI HeKa, INSERM UMR 1138, Sorbonne Université}}

\affil[2]{\orgdiv{Service de Radiologie, Hôpital Armand-Trousseau}, \orgname{APHP}}

\affil[3]{\orgdiv{CMAP, Ecole polytechnique, Institut Polytechnique de Paris}}

\abstract{Estimating accurate high-dimensional transformations remains very challenging, especially in a clinical setting. In this paper, we introduce a multiscale parameterization of deformations to enhance registration and atlas estimation in the Large Deformation Diffeomorphic Metric Mapping framework. Using the Haar wavelet transform, a multiscale representation of the initial velocity fields is computed to optimize transformations in a coarse-to-fine fashion. This additional layer of spatial regularization does not modify the underlying model of deformations. As such, it preserves the original kernel Hilbert space structure of the velocity fields, enabling the algorithm to perform efficient gradient descent. Numerical experiments on \rev{several} datasets, including abnormal fetal brain images, show that compared to the original algorithm, the coarse-to-fine strategy reaches higher performance and yields template images that preserve important details while avoiding unrealistic features. This highly versatile strategy can easily be applied to other mathematical frameworks for almost no additional computational cost.}


\maketitle
\setlength{\headheight}{23pt}%
\addtolength{\topmargin}{-10.60297pt}
\pagestyle{fancy}
\fancyhf{}
\fancyhead{}
\fancyhead[R]{This preprint has not undergone peer review or any post-submission improvements or corrections. The Version of Record of this article is published in \textit{Journal of Mathematical Imaging and Vision}, and is available online at \url{https://doi.org/10.1007/s10851-024-01219-5}}

\section{Introduction} \label{Introduction}

Although the quantitative analysis of anatomical images is an old problem \citep{thompson_1992}, to this day it is still a challenging one, especially given that datasets of clinical images are often small in number and large in dimension. Estimating the transformation that warps an object onto another provides an efficient way of quantifying shape differences, which is the cornerstone idea of Computational Anatomy \citep{Miller2018}. 

\rev{
The choice of the function describing the template-to-subject transformations is of prime importance \citep{Oliveira2014}. To account for the intra- and inter-subject anatomical variability in clinical images, non-linear deformations are mandatory. 
The Large Deformation Diffeomorphic Metric Mapping (LDDMM) setting \citep{Trouve95diffeomorphismsgroups, Miller2002, Christensen1996} is a mathematical framework in which objects are warped through diffeomorphic transformations of the ambient space, i.e. high-dimensional, smooth and invertible functions with smooth inverse that preserve image topology. The group of possible transformations forms a Riemannian manifold of infinite dimension and parameterizes a flexible representation of deformations.

The main advantage of the LDDMM lies in the fact that it constrains diffeomorphisms to be geodesic flows, i.e. the shortest paths between the objects according to a metric regularizing transformations \citep{Trouve95diffeomorphismsgroups}. Therefore, one can compute distances between shapes and perform meaningful statistical analyses on populations of shapes, such as Principal Component Analysis \citep{Vaillant2004}, geodesic regression \citep{Fletcher2011} and atlas estimation \citep{Allassoniere2015}.

In this paper, we are interested in atlas estimation,} a method to model the mean and variability over a collection of images that are instances of the same anatomical object. An estimate of the average shape is given in the form of a template image, which represents the invariants across the population, i.e. shared anatomical features, and the variability is given by deformations from the template space to each subject's space, which express how these common features vary within the population \citep{Grenander1998}.

Atlas estimation has many applications in the field of medical image analysis. The template image can be used as reference to describe average anatomical structures or serve as a tool to automatically segment new subjects. \rev{Within the LDDMM framework, atlases have been used to characterize pathological deviations from normality \citep{gaudfernau2021}, isolate subgroups in a population \citep{Debavelaere2020}, and in a spatio-temporal fashion to characterize pathological changes, such as hippocampal reduction in Alzheimer's disease \citep{Debavelaere2020}.}

\rev{In the LDDMM framework, diffeomorphisms are constructed by integrating time-dependent velocity fields, which results in high computational complexity. Conveniently, the resulting flow of diffeomorphisms is fully determined by the initial state of the system, enabling us to only optimize the initial velocity fields \citep{Miller2006}. However, finding optimal diffeomorphisms involves solving expensive Partial Differential Equations (the “geodesic equations”) \citep{Younes2007} and operating over infinite-dimensional velocity fields. To further improve optimization efficiency, a number of papers introduced sparse finite-dimensional parameterizations of the initial vector fields \citep{durrleman2012, durrleman2011, Sommer2013bis, Zhang2015}. In this paper, we will work within the framework of Durrleman et al. \citep{durrleman2012}, where velocity fields are parameterized by the convolution of momentum vectors attached to control points and regularized by a Gaussian kernel defining a Reproducing Kernel Hilbert Space (RKHS).}

The choice of this regularizer is critical as it restricts the range of transformations defined by the model \citep{Sommer2013}. Specifically, it constrains the deformations occurring on the images to a single scale. A large kernel width is likely to produce smooth but less accurate matches, while a fine kernel will generate more accurate yet unnatural deformations. As clinical images often present high variability at several scales, one might be tempted to increase the number of parameters in the model, i.e. to use many control points and a small kernel. However, fine kernels make large displacements more expensive than small ones, and such over-parameterization will likely trap the optimization procedure in a local minimum, achieving a reasonable numerical solution that is qualitatively bad. 
To overcome such problem, hierarchical algorithms have been widely used in the field of image registration \citep{Oliveira2014, Modersitzki2008}: after solving the registration problem at coarse scales, the solution is transferred to increasingly fine scales to refine the transformation. These strategies avoid more efficiently trapping the algorithm in local minima related to unrealistic transformations. 


\rev{In this paper, we propose to enhance the outcomes of the conventional RKHS-based LDDMM algorithm by employing a coarse-to-fine optimization procedure based on a Haar-like wavelet representation of the initial velocity fields. Our reparameterization can be seen as an additional layer of regularization which does not modify the underlying model of deformation, making our algorithm highly versatile and transferable to other mathematical frameworks. We apply the multiscale strategy to the atlas estimation algorithm of Durrleman et al. \citep{durrleman2012}, which uses a finite parameterization of the velocity fields as a linear combination of RKHS basis elements. This enables us to rely on an efficient numerical scheme to compute the gradients all the while controlling the smoothness of the deformations with our multiscale scheme.}
\rev{We evaluate our algorithm on a registration example and different atlas estimation tasks} and show that our strategy generates more natural template images as well as higher stability regarding the initialization.

This paper is organized as follows. We first explore works related to coarse-to-fine registration and atlas estimation in~\cref{sec:state_of_the_art}, then we recall the LDDMM setting in the case of the finite parameterization of the velocity fields in~\cref{sec:model_of_deformations} and introduce our coarse-to-fine atlas estimation method in~\cref{sec:multiscale_atlas_estimation}; then, we conduct experiments on toy data and fetal brain images in~\cref{sec:experiments} and discuss our results in~\cref{sec:discussion}. 

\section{Related work} 
\label{sec:state_of_the_art}
Even though multiscale image registration has been studied repeatedly in the literature, it has rarely been \rev{evaluated} in the setting of population analysis. As registration is a special case of atlas estimation with a fixed template image, in the following we will review both registration and atlas estimation methods that have a multiscale property.

\rev{
Rooted in the idea of analyzing images across different scales \citep{Witkin1983, Koenderink2004}, hierarchical approaches have played a pivotal role in tasks involving image motion. Traditional optical flow algorithms \citep{Lucas1981, Anandan1989, Enkelmann1998, Terzopoulos1986} were among the first to perform image matching in a coarse-to-fine manner, using pyramids of images of progressively increasing resolution to circumvent the difficulty of capturing large motions and speed up computations. Concurrently, the development of sophisticated tools like Spline pyramids \citep{Unser1993} and wavelet transforms \citep{Mallat1989} provided novel methods for representing images and deformations at different scales. Many popular medical image registration algorithms now use multiresolution representations of deformations, most of which are combined with Gaussian smoothing of the images: the Demons algorithm \citep{Thirion1998} applies Gaussian filters to the vector fields, IRTK \citep{Rueckert1999, Schnabel2001} parameterized B-splines with grids of increasing resolution, ROMEO \citep{Hellier2007} uses an adaptive multigrid strategy to estimate optical flow, DARTEL \citep{ASHBURNER2007} also uses a multigrid method to estimate single-flow velocity fields and Syn \citep{avants2008} estimates symmetric diffeomorphisms with Gaussian smoothing of velocity fields and later on B-spline regularization \citep{Tustison2013}.}

In the following, we will only review multiscale strategies related to our method, i.e. that are multiscale with regard to the deformation field.

\subsection{Multiscale representations of deformations}

Signal representation provides numerous mathematical tools to represent signals, and several registration algorithms modelled deformations using basis functions from the Fourier and Wavelet transforms. \rev{These multiscale representations typically aim to improve computational efficiency by playing on the relation between a representation and a regularization, which allows using lower dimensional parameterization of the velocity fields along the trajectory. Examples of this approach made use of the Fourier transform \citep{Ashburner1999, Zhang2015, Zhang2017}, Galerkin method \citep{Mang2015}, or orthogonal decomposition \citep{Zhang2019}. Our work falls into a different category, which aims to modify an initial method marginally, to reach better local minima without changing the considered representation or regularization of the LDDMM. The solution is to construct a sequence of nested problems: in our case, a multiscale representation of only the initial velocity field. The inherent scale information within the new parameters is leveraged to solve the optimization problem in a coarse-to-fine manner. This idea is already present in Christensen et al. \citep{Christensen2001} who parameterized displacement fields by a Fourier basis with an adapted regularization term and estimated frequency coefficients in a coarse-to-fine fashion.} However, experiments showed that modelling deformations using a wavelet basis provides better spatial regularization compared to using a Fourier basis \citep{Amit1994}. This advantage stems from the fact that wavelets not only capture the frequency of the signal but also its location and orientation, making them a more natural fit for hierarchical optimization strategies.

Wavelet-based deformation models cover various wavelet types (e.g. Haar, Cai Wang and $(BV,L^2)$) and deformations (e.g. displacement vectors, B-splines and elastic deformations). 
Displacement vectors \citep{Wu2000, Cathier2006} and later free form deformations \citep{Sun2014} were described and optimized in a multi-resolution fashion through the Cai Wang wavelet \citep{Cai1996}. Gefen et al. \citep{Gefen2004} modeled elastic deformations using finite-supported, semi-orthogonal wavelet functions. Topology-preserving displacement fields were modelled by polynomial spline basis functions and controlling the Jacobian of the transformation \citep{Musse2001, Noblet2005}. Recently, wavelets $(BV,L^2)$ were employed to generate a hierarchical representation of diffeomorphisms in the hyperelasticity framework and perform coarse-to-fine optimization \citep{debroux2021}.

\subsection{Multiscale registration in the LDDMM framework} 
\label{sec:state_of_the_art_lddmm}
In the LDDMM framework, the choice of the spatial regularizer restricts the range of possible deformations to those occurring at a single scale, which often proves unrealistic \citep{Risser2011}. Thus, a variety of papers have focused on increasing the flexibility of the deformation model. These strategies can be broadly classified into two categories:  the first one \textit{simultaneously} estimates coexisting flows of different scales \citep{Risser2011, Bruveris2011, Sommer2013, Gris2018, Tan2016}, and the second one composes multiple scale flows which are estimated \textit{sequentially} \citep{Modin2018, Miller2020}. 

Risser et al. \citep{Risser2011, Bruveris2011} first introduced a multi-kernel extension of the LDDMM framework in which the deformation flow is defined by a weighted sum of Gaussian kernels whose widths are specified by the user. Weights are tuned in a semi-automatic manner during a pre-registration step.  In this framework, the RKHS structure of the velocity fields is lost, and a new definition of the norm is used to ensure efficient computation of the flow. 
A spatially-varying version of this framework, the kernel bundle, used sparsity priors to allow the weights of the kernel mixture to vary across spatial locations \citep{Sommer2013}. Even though this algorithm proved efficient on the registration of landmark points, the increase in computational cost restricts its application to registration problems involving few parameters. Using an algorithm for the multiresolution decomposition of surfaces, the kernel bundle framework was also combined with a coarse-to-fine strategy wherein the resolution of cortical surface meshes is progressively increased along with that of the deformation field \citep{Tan2016}. 
Multi-kernel approaches were further combined with deep learning optimization in order to learn a local regularizer from the data \citep{Niethammer2019, Shen2019}. However, these methods increase significantly the complexity of the mathematical model, and several optimization procedures are required to tune the networks parameters, the kernel pre-weights and the deformation parameters.

%
A related approach, based on modular deformations, enables the user to impose spatially-varying constraints on the deformation field \citep{Gris2018}. Diffeomorphisms are built by superimposing deformations modules that  encode local geometrical transformations.
As in the kernel bundle framework, the space of vector fields is equipped with an adapted norm. The need for prior knowledge about the deformation modules limits the practical application of the algorithm.

A second and less explored axis of research constructed a hierarchical representation of deformations, based on non-coexisting vector flows of increasing resolution, which are estimated independently and then composed.
In a theoretical paper, Modin et al. \citep{Modin2018} extended $(BV,L^2)$ wavelets to express diffeomorphisms as a composition of deformations of increasingly fine scales, which can be seen as a series of LDDMM steps. Despite the potential of this approach, the authors did not perform numerical experiments. 
A similar approach \citep{Miller2020} constructed diffeomorphisms by composing a series of multiscale vector fields, which enables to progressively refine the deformation. Contrary to multi-kernel approaches, such strategies perform optimization in successive RKHS of increasingly finer resolution, in the spirit of coarse-to-fine strategies. 
\\

As we shall see in the following, our coarse-to-fine approach is more closely related to the one that composes multiple scale flows, in the sense that we perform optimization sequentially in sub-spaces of increasing resolution. 
However, our algorithm differs from the previous ones by the fact that our work brings changes to the optimization procedure rather than the deformation model: the multiscale structure is only used for the initial velocity field and so that the velocity fields are still defined, at core, by a single-scale RKHS. This simplifies the implementation of our algorithm while preserving the efficient optimization scheme of Durrleman et al. \citep{durrleman2012}.


\section{Model of diffeomorphic deformations}\label{sec:model_of_deformations}

\subsection{Large Deformation Diffeomorphic Metric Mapping} 
\label{sub:large_deformation_diffeomorphic_metric_mapping_lddmm_}

In the following, we consider a set of $N$ images $(I_i)_{ 1\leq i  \leq N}$ of dimension $d$. We assume that each image $I_i$ is a smooth deformation of a template $I_{ref}$ plus an additive random white noise ${\epsilon_i}$:
\begin{eqnarray}
\label{eq:model}
I_i = I_{ref} \circ \flowphi_i^{-1} + \epsilon_i,  \ \forall i \in [1, n]
\end{eqnarray}
where ${\flowphi_i}$ is the $i^{th}$ template-to-subject deformation, and ${I_{ref} \circ \flowphi_i^{-1}}$ denotes the action of the diffeomorphic deformation on the template.

In atlas estimation, one seeks to estimate the template image $I_{ref}$ and the $N$ template-to-subject deformations ${(\flowphi_i)_{ 1\leq i  \leq N}}$. Note that registration is a specific case of atlas estimation where the template image is fixed and $N=1$.

We choose to work in the LDDMM setting \citep{Trouve95diffeomorphismsgroups, Miller2002}, in which objects are deformed via deformations of the whole ambient space. This framework generalizes the linearized deformation setting in order to define diffeomorphic deformations that are invertible and smooth. Diffeomorphisms are constructed by integrating linearized deformations over time, which are considered as infinitesimal steps, according to the differential flow equation: 
\begin{equation}
	\label{eq:flow_equation}
	\left\{
	\begin{aligned}
	\frac{d\locpix(t)}{dt} &= v_t(\locpix (t))\\
	\locpix(0) &=  \locpix_0\,.
	\end{aligned}\right.
\end{equation}
where $v_t \in V$ is an instantaneous velocity field belonging to a Hilbert space $V$ and $x$ can be seen as a particle moving along the curve $\locpix(t)$ in the domain of interest $D$.

This model builds a flow of diffeomorphisms 
$\flowphi_t: \locpix_0 \To \locpix(t)$ $\forall t\in[0,1]$. 
The diffeomorphism of interest $\flowphi_1$ is defined as the end point of the path $x(t)$, i.e.:
\begin{equation*}
	\forall \locpix_0\in D, \ \ \flowphi_1( \locpix_0) = \locpix(1)\,.
\end{equation*}
Note that for any time $t$, $\flowphi_t$ is indeed a diffeomorphism  provided that the velocity field is regular enough, i.e. that it is continuous squared integrable. 

\subsection{Discrete parameterization of diffeomorphisms}

Finally, we need to define an appropriate norm $||.||_V$ for the Hilbert space $V$. To this end, we restrict ourselves to vector fields that belong to a RKHS \citep{Berlinet2004} $V$ defined by a kernel $K_g$. 
We also rely on the work of Durrleman et al. \citep{durrleman2012} to introduce a discrete parameterization of the velocity fields: we assume that the initial velocity field $v_{0}$ can be decomposed as a \textit{finite} linear combination of the RKHS basis vector fields. The weights of the decomposition of a given deformation onto this basis are given by a set of momentum vectors $(\alpha_{k}(0))_k$ attached to $k_g$ control points $(c_{k}(0))_k$.
\begin{align*}
    v_{0,i}(\locpix) =  \sum\limits_{k=1}^{k_g} \Kg (\locpix,c_{k,i}(0)) \, \alpha_{k,i}(0).
\end{align*}
In this work, ${\Kg}$ is the Gaussian kernel:
$\Kg(x, y) = \exp(\frac{- \|x - y\|^2}{\sigma_g}^2) I_d$, 
with $\sigma_g$ the kernel width and $I_d$ the identity matrix.

 The structure of the RKHS is such that the vector fields are continuous and squared integrable. As we shall see in \Cref{sec:optimization}, this will turned out to be important for performing optimization in a finite-dimensional setting.

It has been proved \citep{Miller2006} that the vector fields that define geodesic deformations with respect to the norm $\int_{0}^1 \|v_{t,i}\|_{V}^2 \, dt$ keep the same structure along time and write according to:
\begin{equation}
	\label{eq:velocity}
	v_{t,i}(\locpix) =  \sum\limits_{k=1}^{k_g} \Kg (\locpix,c_{k,i}(t)) \alpha_{k,i}(t)\,,
\end{equation}
where for any time $t$, $\alpha_{k,i}(t)$ is the $k^{th}$ momentum vector related to subject $i$ and attached to the control point $c_{k,i}(t)$.

 Furthermore, the trajectory of the control points  $(c_{k,i}(t))_k$ and momentum vectors $(\alpha_{k,i}(t))_k$ is described by the Hamiltonian system equations \citep{Miller2006}:
 \begin{equation}\label{eq:HamSyst}
	\left\{
	\begin{aligned}
		\frac{d\gk(t)}{dt} & = \sum_{l=1}^{k_g} \Kg (\gk(t),c_l(t))\alpha_l(t)\\
		\frac{d\alpha_k(t)}{dt} &= -\sum_{l=1}^{k_g} d_{c_k(t)} (\Kg (\gk(t),c_l(t))\alpha_l(t)^t\alpha_k(t)\
	\end{aligned}
	\right.
\end{equation}
with initial conditions $\gk(0)$ and $\alpha_k(0)$ for all $1\leq k\leq k_g$.  

 Finally, one verifies that the kinetic energy along geodesic paths is preserved over time, i.e. $\forall t\in[0,1], \|v_{t,i}\|_V= \|v_{0,i}\|_V\nonumber$.
This implies that a geodesic transformation is fully parameterized by the initial velocity field $v_{0,i}$. Hence, estimation of the diffeomorphism $\flowphi_{1,i}$ boils down to a geodesic shooting problem. The system is deterministic, and we only need to optimize the initial conditions $\alpha_{0,i}=(\alpha_{k,i}(0))_{k}$ and $c_{0,i}=(c_{k,i}(0))_{k}$ for each subject $i$ along with the template image $I_{ref}$.

\subsection{Optimization}\label{sec:optimization}

\begin{figure*}[t]
\centering
    \includegraphics[width=0.95\linewidth]{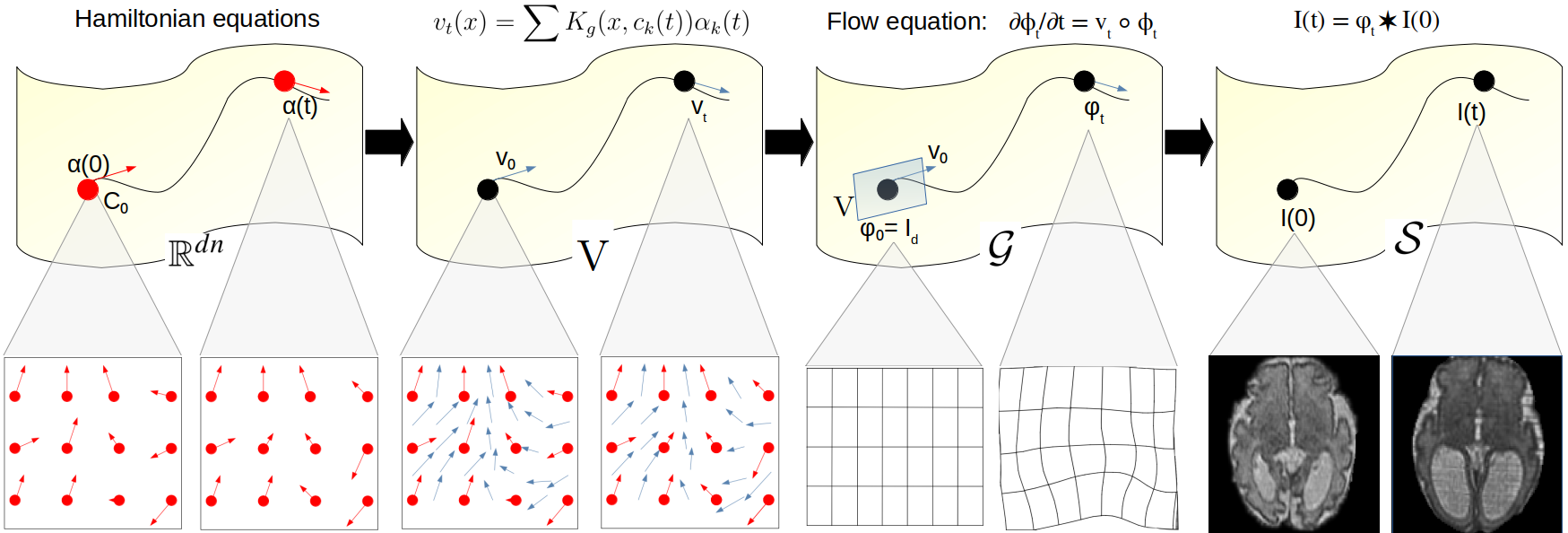}
    \caption{Computing shape deformations from an initial image $I(0)$ and a set of initial momentum vectors $\alpha(0)$. 1- Given $\alpha(0)$ and the initial control points $c_0$, integration of the Hamiltonian (\Cref{eq:HamSyst}) gives the trajectory of the momentum vectors $\alpha(t)$. 2- The velocity field $v_t$ is computed by interpolating the momentum vectors with \Cref{eq:velocity}. 3- Integration of the flow equation (\Cref{eq:flow_equation}) gives a flow of diffeomorphisms $(\flowphi_t)_{t \in [0,1]}$. 4- Finally, $\flowphi_t$ is applied to the object $I(0)$, giving the deformed image $I(t)$ at any time $t$.}
    \label{fig:lddmm_optim}
\end{figure*}


In this work, the position of the control points is fixed on a regular grid of spacing $\sigma_g$. 
We seek to minimize a cost function expressing a trade-off between matching accuracy and regularity of the transformation(s). With the parameterization introduced in the previous section, the cost function $E$ is a function of $I_{ref}$ and $(\alpha_{0,i})_i$:
\begin{align}
	\label{eq:objectiveparam}
  &E(I_{ref},(\alpha_{0,i})_{1\leq i\leq N})= \quad \nonumber \\
 & \sum_{i=1}^{N} \left( \frac{d(I_i, I_{ref}  \circ \flowphi_{1,i}^{-1})^2}{\sigma_\epsilon^2}  + \|v_{0,i}\|_V^2 \right),
\end{align}
The regularity term is the total kinetic energy along the geodesic path related to $\flowphi_i$. 
With the discrete parameterization chosen for the velocity fields (i.e. \Cref{eq:velocity}), the norm $\|.\|_V$ can easily be computed and thus
\begin{equation}
\label{eq:total_kinetic_energy}
 \|v_{0,i}\|_V^2 = \sum\limits_{j=1}^{k_g} \sum\limits_{k=1}^{k_g} \alpha_{j,i}(0) \, \Kg (c_{j,i},c_{k,i}) \, \alpha_{k,i}(0)\, .
\end{equation}
\rev{Optimization is performed through gradient descent.} Thus, we need to compute the gradient with respect to all the parameters.
This is not a trivial task, but an efficient numerical scheme has been proposed in Durrleman et al. \citep{durrleman2012}. 
This algorithm relies heavily on the fact that the norm of the vector fields is the one of the RKHS and that the vector fields that are
solutions to the problem remain a \textit{finite} sum of kernels at all times (\Cref{eq:velocity}). This enables us to solve a finite-dimensional problem even though the functions we are looking for have infinite dimension.   
Given a set of momentum vectors $\alpha_0$, this strong structural property
allows to efficiently deform shapes and compute the cost function by performing the following steps, which are illustrated in \Cref{fig:lddmm_optim}: integrating \Cref{eq:HamSyst} gives the evolution of the momenta over time, the velocity field $v_t$ at any time $t$ is computed with \Cref{eq:velocity}, the flow of diffeomorphisms $(\flowphi_t)_{t\in[0,1]}$ is obtained by solving \Cref{eq:flow_equation}, and the template image $I_{ref}$ is deformed with the flow: $\flowphi_1 \star I_{ref}$. It is then straightforward to infer the distance between the deformed template and the target objects and to compute the regularity term with \Cref{eq:total_kinetic_energy}.
Finally, one can compute the gradient of the cost function with respect to the template, and use a backward integration along time to compute the gradient with respect to the momentum vectors.
These steps correspond to the lines \ref{line:gradients_1}-\ref{line:gradients_2} in  \cref{alg:algorithm_ctf}.

The algorithm is publicly available as part of the open-source software Deformetrica \citep{Deformetricasoftware}.

\subsection{Unrealistic local minima issue}



There is a dependency between the scale of the kernel and the number of parameters, as a constant vector field has to be well approximated by the finite sum: this imposes that the scale $\sigma_g$ should be related to the distance between control points. Here, the kernel $K_g$ not only controls the \textit{regularity} of the deformations, but also the \textit{number of parameters} to optimize.
This leaves the user to find the balance between a large kernel, (i.e. few control points), which generates smooth but less accurate transformations, and a smaller one, which penalizes large displacements and increases the risk of converging towards unrealistic solutions. 

Further, as the optimization problem is not convex, the gradient descent algorithm converges towards a solution that depends on the initialization. The more the numbers of subjects and parameters are high, the more complex the energy landscape of the problem becomes.
In the original approach of Durrleman et al. \citep{durrleman2012}, a first step towards the estimation of multiscale deformations was taken by optimizing both the \textit{number} and \textit{position} of the control points. In areas with low variability, points are inactivated by a L1 sparsity prior on the momenta. Unfortunately, the width of the kernel remains the same and in practice, the vector fields can only be set to zero in image areas of null intensity. Another issue is that the position of the control points cannot change significantly if the distance between points is low.
Yet, it is exactly when the number of parameters is high that the risk of converging towards an unrealistic local minimum is higher. 


One idea would be to change locally the scale of the kernel so that we can estimate non evenly smooth vector fields, i.e. vector fields with spatially varying scales.
Even if such parameterization can be written, one looses the RKHS structure and thereby the ease of computation.

In this paper, we wish to address these two related issues: the dependency of the algorithm on the initialization, which restrains the number of parameters that can be properly optimized, and the difficulty in estimating vector fields that have \textit{locally} varying regularity.
In the next section, we will describe a reparameterization of the vector fields which enables us to impose smoothness constraints on the deformations and progressively relax them in a coarse-to-fine fashion. In this way, the algorithm can cope with non evenly smooth transformations while using a small kernel and remaining in the original RKHS setting. 

\section{Multiscale deformations}
\label{sec:multiscale_atlas_estimation}

In this section, we propose a multiscale optimization procedure based on a Haar-like wavelet representation of the initial velocity fields. This strategy has the advantage of making the algorithm less dependent on the initialization while favoring more multiscale  deformations.
We rely on the finite parameterization of the velocity fields as a linear combination of RKHS basis elements \citep{durrleman2012} within the LDDMM setting. Importantly, our strategy enables us to preserve this structural assumption and the efficient numerical scheme that follows. We will show that our algorithm generates more natural template images as well as higher stability regarding the initialization.

\subsection{Overview}

\begin{figure*}[h]
    \centering
    \subfloat[Classical LDDMM]{\includegraphics[width=0.2395\textwidth]{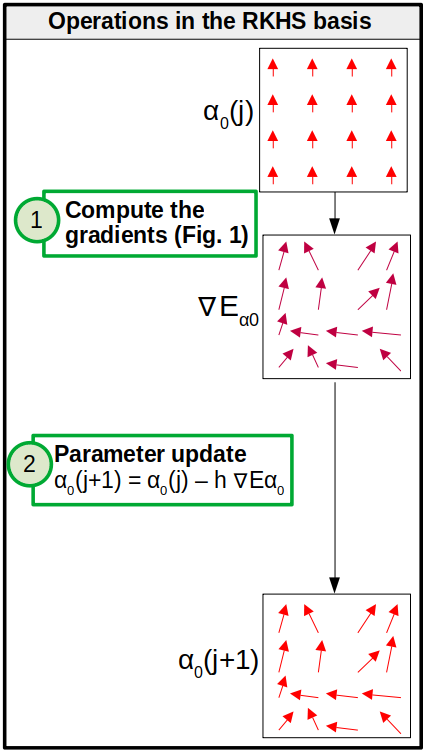}} \hspace{1.5cm}\subfloat[LDDMM $+$ multiscale optimization]{\includegraphics[width=0.63\textwidth]{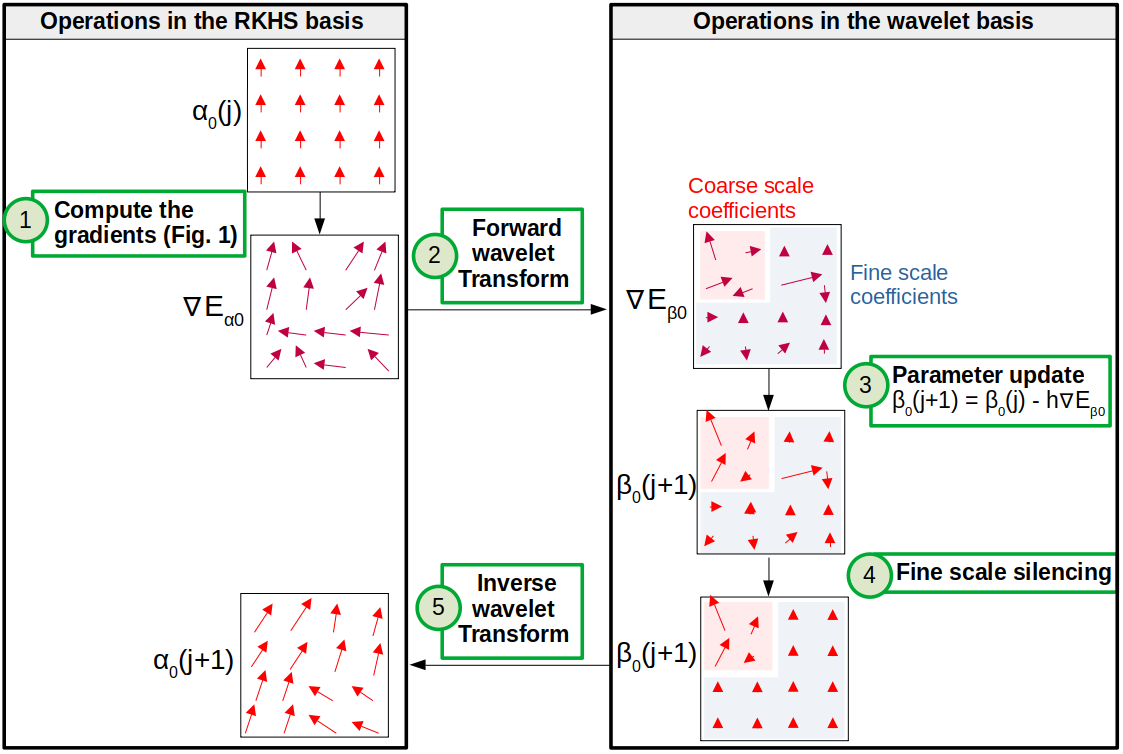}}
    
    \caption{Overview of the original LDDMM algorithm \citep{durrleman2012} and the coarse-to-fine optimization strategy in dimension 2. For the sake of clarity, only a \textbf{single} gradient descent iteration is presented in each panel. The original algorithm iterates between two classical steps: computing the gradients (1) and updating the parameters (2). In the multiscale strategy, step (1) is preserved and step (2) is replaced by another procedure involving updating the representation of the momentum vectors in the wavelet basis. Red arrows denote the coordinates of the momentum vectors either in the RKHS basis (for $\alpha_0(j)$) either in the wavelet basis (for $\beta_0(j)$). Purple arrows denote the gradients of the cost function $\nabla E_{\alpha_0}$ and $\nabla E_{\beta_0}$ with regard to $\alpha_0$ and $\beta_0$ (respectively).}
    \label{fig:ctf_overview}
\end{figure*}

For the sake of clarity, we provide the reader with an overview of our multiscale strategy.
In the \textit{original} algorithm for atlas estimation, we optimize a cost function $E$. Two types of parameters are optimized through gradient descent: the template image $I_{ref}$ and $N$ sets of momentum vectors $\alpha_{0,i}$ that parameterize the template-to-subjects velocity fields $v_{0,i}$. 
The \textit{original} optimization, summarized in \Cref{fig:ctf_overview}, iterates between two classical steps (note that the subscripts $i$ have been dropped for clarity):

\begin{enumerate}
\item Computation of the gradients $\nabla E_{\alpha_0}$ and $\nabla E_{I_{ref}}$ (\Cref{sec:optimization})
\item Parameters update: $\alpha_0(j+1) \leftarrow \alpha_0(j) - h \times \nabla E_{\alpha_0}$ and $I_{ref}(j+1) \leftarrow I_{ref}(j) - h \times \nabla E_{I_{ref}}$, with $h$ the step size.
\end{enumerate}

In the \textit{multiscale} strategy, summarized in \Cref{fig:ctf_overview}, we only modify the parameters update step of the momentum vectors by replacing it with the following steps:

\begin{enumerate}[label=(\roman*)]
\item We use the wavelet transform to obtain a multiscale representation of the gradient of $E$ with regard to the momenta: $\nabla E_{\beta_0} \leftarrow FWT(\nabla E_{\alpha_0})$. (This reparameterization is detailed in \Cref{sec:reparameterization})
\item The coordinates of the initial velocity fields $v_0$ in the wavelet basis are updated: $\beta_0(j+1) \leftarrow \beta_0(j) - h \times \nabla E_{\beta_0}$.
\item The wavelet coefficients in $\beta_0(j+1)$ whose scale is smaller than a current scale $S_j$ are set to zero.
\item The coordinates of $v_0$ in the RKHS basis are recovered with the Inverse Wavelet Transform: $\alpha_0(j+1) \leftarrow IWT(\beta_0(j+1))$; $I_{ref}$ is updated as in the original step (2). (The coarse-to-fine optimization steps are detailed in \Cref{ref:coarse_to_fine_atlas_estim}.)
\end{enumerate}


\subsection{Reparameterization of the initial velocity fields}\label{sec:reparameterization}

In this section, we introduce a multiscale reparameterization of the initial velocity fields based on the Haar wavelet transform. The choice of the Haar wavelet is motivated by its ease of implementation and the orthogonality of its transform. We introduce the definition and properties of the continuous Haar Wavelet representation, recall how this construction can be extended to a representation of discrete signals defined on a grid, and we demonstrate how this can be used to obtain a Haar-like representation of the initial velocity fields. 

\subsubsection{The continuous Haar wavelet transform}

\begin{figure}[ht]
\centering
\includegraphics[width=0.98\linewidth]{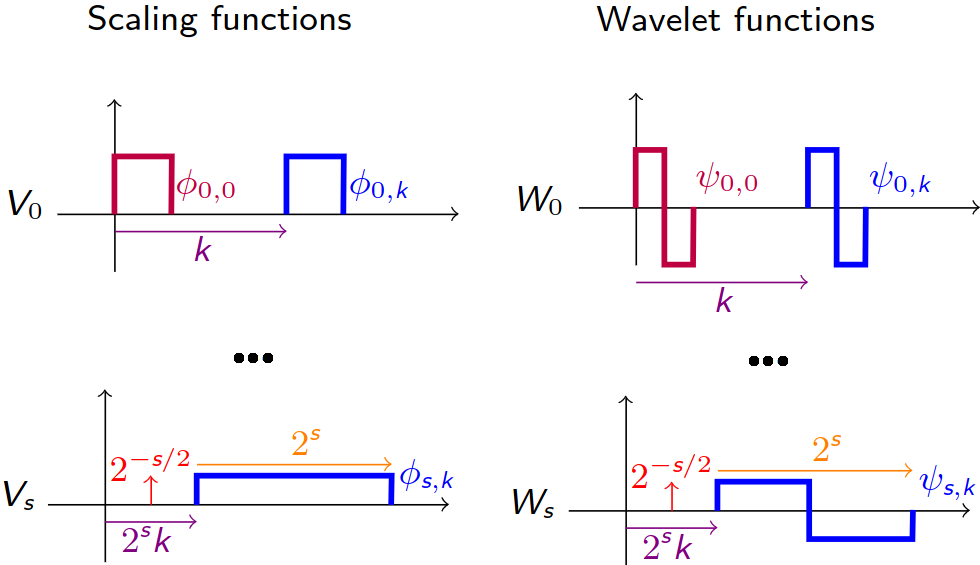}
\caption{The 1D Haar wavelet. At scale 0, $V_0$ is the space of piecewise constant functions of size $1$. A basis for $V_0$, i.e. $\left\{\flowphi_{0,k} \right\}_k$, is obtained by translating the scaling function $\flowphi_{0,0}=\psi^L$ by factors $k \in \Z$.  A basis for $W_0$, i.e. $\left\{\psi_{0,k} \right\}_k$, is obtained by translating the mother wavelet function $\psi_{0,0}=\psi^H$ by factors $k \in \Z$. At scale $s$, a basis for the space $V_s$, i.e. $\left\{\flowphi_{s,k} \right\}_k$,  is obtained by dilating $\flowphi_{0,0}$ by $2^s$, translating it by $2^sk$ and normalizing it by $2^{-s/2}$. A basis for $W_s$, i.e. $\left\{\psi_{s,k} \right\}_k$, is obtained by performing the same operations on $\psi_{0,0}$.}
\label{fig:wavelets}
\end{figure}

Here, we describe the decomposition of a real signal $f$ defined on the $d$-dimensional space $\mathbb{R}^d$ into a Haar Wavelet basis \citep{Mallat1989, Mallat2009}.
The wavelet representation decomposes $f$ into a linear combination of basis functions which have different \textit{resolutions}, \textit{locations} and \textit{orientations}. This representation relies on a collection of embedded spaces $V_s$ that contain functions said of scale $s$. 

In the case of the Haar Wavelet, $V_s$ is the space of piecewise constant functions on a regular grid of size $2^s$. Any function $f$ can be approximated in this space by computing a local average: the mean value in each sub-square of the grid.

We define the $d$-dimensional \textit{scaling function} $\phi$ by \[\phi(x) = \prod_{i=1}^d \psi^L(x_i)\]
where $x \in \mathbb{R}^{d}$ and $\psi^L$ is the 1D piecewise constant function
\[\psi^L(z) =
    \begin{cases}
      1 & \text{for ${0 \leq z < 1}$}\\
      0 & \text{otherwise}
    \end{cases}. \]
As illustrated in \Cref{fig:wavelets}, approximating $f$ at scale $s$ amounts to projecting $f$ onto the space spanned by the orthonormal family
\[\left\{ \phi_{s,k}(x)= 2^{-sd/2} \phi(2^{-s}x-k) \right\}_{k \in \mathbb{Z}} \]
where $\phi_{s,k}$ is the scaling function rescaled by $2^s$ and then translated by $k 2^s$. The factor $2^{-sd/2}$ is a normalization factor that ensures that the $\phi_{s,k}$ have unit energy. 

When transitioning from the approximation at scale $s$ to the approximation at the coarser scale $s+1$, some details of $f$ are lost. These details belong to the orthogonal complement $W_{s+1}$ of the space $V_{s+1}$ in $V_{s}$. 
A basis of this space can be obtained by defining the $d$-dimensional oriented wavelet functions
\[\psi^{o}(x) = \prod_{i=1}^d \psi^{o_i}(x_i)\]
where $x \in \mathbb{R}^{d}$, $o \in \{H,L\}^d$, $\exists \: i, o_i = H$ and
\[\psi^H(z) =
    \begin{cases}
      1 & \text{for ${0 \leq z < 0.5}$}\\
      -1 & \text{for ${0.5 \leq z < 1}$}\\
      0 & \text{otherwise}
    \end{cases}.\] 
As illustrated in \Cref{fig:wavelets}, an orthonormal basis of $W_{s}$ is given by
\[
\left\{ \psi_{s,k}^{o}(x)= 2^{-sd/2} \psi^{o}(2^{-s}x-k) \right\}_{k \in \mathbb{Z},\: \exists i\:o_i = H}.
\]
where $\psi^{o}_{s,k}$ is the wavelet function of orientation $o$ rescaled by $2^s$ and then translated by $k 2^s$. 
Note that functions $\psi^L$ and $\psi^H$ act respectively as low and high pass filters. Their combination yields oriented high pass filters, e.g., for $d=2$, there exist three wavelet functions $\psi^{HL}$, $\psi^{LH}$ and $\psi^{HH}$, that express details of the signal along vertical, horizontal and diagonal orientations (respectively).

We can thus decompose any function $f$ in $V_{s}$ in the two following ways:
\begin{align*}
f &= \sum_{k}^{} a_{s,k} \; \phi_{s,k}\\
& = \sum_{k}^{} a_{s+1,k}\;  \phi_{s+1,k} + \sum_{o,k}^{} d^{o}_{s+1,k}\;  \psi^{o}_{s+1,k}
\end{align*}
As $\sum_{k}^{} a_{s+1,k} \phi_{s+1,k}$ belongs by construction to the space $V_{s+1}$, we can further decompose it into
a projection onto $V_{s+2}$ and a projection onto $W_{s+2}$. Repeating this scheme up to scale $S$ leads to the following multiscale decomposition of $f$:
\[
f= \sum \limits_{k}^{} a_{S,k}\;  \phi_{S,k} + \sum \limits_{s'=s+1}^{S} \sum \limits_{o,k} d^{o}_{s',k} \;  \psi^{o}_{s',k}.
\]
The classical wavelet construction is concluded by letting $s$ go to $-\infty$, enabling one to decompose any measurable bounded function in such bases. 

More importantly, going from the decomposition in $V_s$ to the decomposition in the spaces $V_S$ and $(W_{s'})_{s<s'\leq S}$ corresponds to a change of basis and thus to a discrete operation going from the coefficients $(a_{s,k})_k$ to the coefficients $(a_{S,k})_k$ and $(d^o_{s',k})_{s<s'\leq S, k,o}$. This transformation is called the Forward wavelet transform (FWT) and its inverse the Inverse wavelet transform (IWT). Both can be computed directly on the coefficients without relying on the continuous basis functions.

\subsubsection{Haar wavelet applied to grids}

\begin{figure*}[tb]
\subfloat[Classical decomposition into unit functions]{\includegraphics[width=0.9\linewidth]{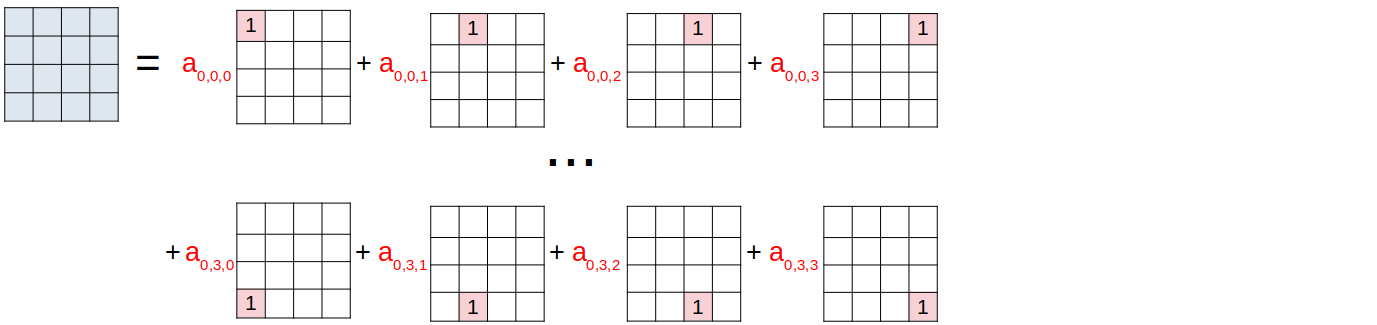}} \\
\subfloat[Multiscale Haar wavelet basis]{\includegraphics[width=0.9\linewidth]{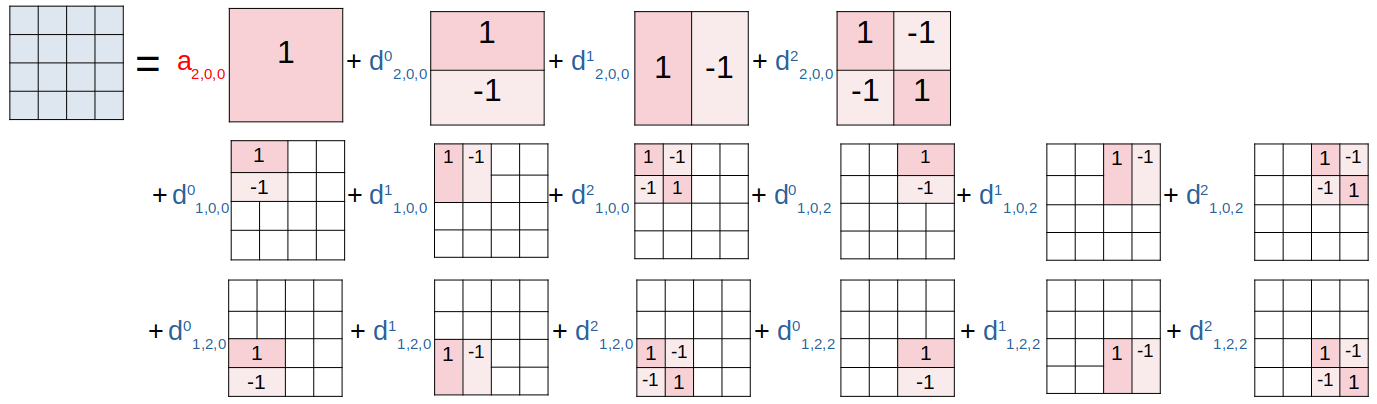}}     
\caption{Decomposition of a 4-by-4 grid in two bases. The letters a and d refer to approximation and detail coefficients, respectively. Subscripts indicate the \textit{scale} of the coefficient and the x-y \textit{position} of the related wavelet function. Superscripts indicate the \textit{orientation} of the wavelet function. In the grids, empty cells denote null values.}
\label{fig:haar_on_vectors}
\end{figure*}

In our algorithm, rather than using the continuous Haar wavelet decomposition, we have implemented the related discrete Haar wavelet decomposition on a $d$-dimensional grid
$[\![0,K_1]\!] \times ... \times [\![0,K_d]\!]$, illustrated in \Cref{fig:haar_on_vectors} in 2D. 
This transform corresponds to a discrete change of basis related to the continuous Haar transform.

More precisely, to any discrete function $(a_{0,k})_k$ on the grid, 
one can associate the continuous function $f$ of $V_0$  defined by
\[ f = \sum_{k \in [\![0,K_1]\!] \times ... \times [\![0,K_d]\!]} a_{0,k} \;  \phi_{0,k},
\]
where the $(a_{0,k})_k$ are interpreted as the approximation coefficients at scale $0$ of function $f$. 

By construction, $f$ is a piecewise constant function on the related continuous grid.
Using the FWT algorithm up to scale $S$, this function can be decomposed as follows:
\rev{
\begin{align*}
    f = \sum \limits_{k}^{} a_{S,k}\;  \phi_{S,k} + \sum \limits_{s=1}^{S} \sum \limits_{o,k} d^{o}_{s,k} \;  \psi^{o}_{s,k}.
\end{align*}}
One verifies that $a_{S,k}=0$ when $k_i<0$ or $k_i2^S>K_i$ and $d_{s,k}=0$ when $k_i<0$ or $k_i2^{s}>K_i$. Thus, these sums have only a finite number of coefficients.
Further, when the $K_i$ are powers of 2, i.e. $K_i=2^{\zeta_i}$, we impose that the decomposition cannot exceed a maximum scale $\smax = \min_i (\zeta_i)$, so that the scaling function support remains within $\prod_i [0,K_i]$. Thus, when $S \leq \smax$, the previous equation reduces to:
\rev{
\begin{align}
    f = \sum \limits_{0 \leq k_i < 2^{\zeta_i}-S}^{} a_{S,k}\;  \phi_{S,k} \nonumber\\
    + \sum \limits_{s=1}^{S} \sum \limits_{o, \, 0 \leq k_i < 2^{\zeta_i-s}} d^{o}_{s,k}\;  \psi^{o}_{s,k},
    \label{eq:haar_decomposition_grid}
\end{align}}
which corresponds exactly to an orthonormal change of basis.

To summarize, from the coefficients $(a_{0,k})_k$, one can compute the coefficients $(a_{S,k})_{0 \leq k < 2^{\smax-S}}$ and $(d^{o}_{s,k})_{1 \leq s \leq S,0 < k < 2^{\smax-s},o}$ with the FWT algorithm and perform the inverse operation with the IWT algorithm.

We denote by $M_{FWT}$ and $M_{IWT}$ the $\prod_i K_i \times \prod_i K_i$ matrices associated to the transformations FWT and IWT. $M_{FWT}$ and $M_{IWT}$ are related to one another: we have obviously $M_{FWT} * M_{IWT}$ = $Id$ and, as these are orthonormal transforms, we also verify that $M_{IWT}^T$ $=$ $M_{FWT}$. Furthermore, there exist fast implementations of these transforms that have a linear complexity with respect to the number of coefficients \citep{Mallat2009}.


In this work, we have implemented the FWT using the fast lifting scheme described in \citep{Mallat1989} (Chapter 7.8). This scheme is strictly equivalent to the previous description when the grid size is a power of $2$, with the difference that it handles non-dyadic grids through improved computations at the boundaries and that FWT and IWT remain orthogonal transforms. 
The pseudocode for the algorithms FWT and IWT can be found in Appendix (Algorithms \ref{alg:FWT} and \ref{alg:IWT}, respectively), along with a brief explanation and an illustration of how the FWT algorithm operates on a non-dyadic grid (\Cref{fig:algorithm_wavelet}).

Note that these algorithms are available in a public Git repository\footnote[1]{\url{https://github.com/fleurgaudfernau/Deformetrica_multiscale/}}.




\subsubsection{Preservation of the RKHS structure of the velocity fields}

The previous subsection transposed the classical Haar description of continuous functions to discrete functions defined on a grid. To implement a coarse-to-fine initialization approach for atlas estimation, we will describe the initial momentum vectors not through their values on the grid but by the decomposition of these values in the discrete Haar basis. The multiscale structure will be used to obtain a smooth initial field by setting fine-scale coefficients to $0$. As we will show in this section, this does not change the fact that the initial vector fields $v_{0,i}$, and thus all the vector fields $v_{t,i}$, are by construction finite combinations of the RKHS kernel $K_g$. 

If we consider two spaces $V_s$ and $V_{s'}$ with $s \geq s'$, we observe that any function $\phi_{s,k} \in V_{s}$ or $\psi^{o}_{s,k} \in V_{s}$ is also in the finer space $V_{s'}$. In particular, if we take $s'=0$,  $\phi_{s,k}$ and $\psi^{o}_{s,k}$ can be decomposed as a linear combination of $\phi_{0,k'}$:
\rev{\begin{align}
    \phi_{s,k} &= \sum_{k'} \gamma_{s,k,k'} \; \phi_{0,k'}\label{eq:phidec}\\
    \psi^{o}_{s,k} & = \sum_{k'} \gamma^{o}_{s,k,k'}\;  \phi_{0,k'}\label{eq:psidec}
    \end{align}}
    where $\gamma_{s,k,k'}$ and $\gamma^{o}_{s,k,k'}$ are some fixed real numbers.
    
\rev{    Thus, for a function $v=\sum_{k} a_{0,k} \; \phi_{0,k}$, its wavelet coefficients $((a_{S,k})_{k},(d^{o}_{s,k})_{1 \leq s \leq S, o, k})$ can be computed as
    \begin{align*}
    a_{S,k} &= \sum_{k'} \gamma_{S,k,k'} \; a_{0,k'} \\
    d^{o}_{s,k} & = \sum_{k'} \gamma^{o}_{s,k,k'}\; a_{0,k'}  
\end{align*}}

Further, we recall that we optimize shape transformations by using the parameterization of Durrleman et al. \citep{durrleman2012}, i.e.
the initial geodesic vector field of a given subject $i$ is defined as a finite linear combination of identical Gaussian kernels that are evaluated at an initial set of points $(c_{k,i}(0))_k$. 
In our scheme, we set the initial points on a grid so that this initial velocity field writes as
\begin{align}
    v_{0,i}(\locpix) =  \sum\limits_{k \in [\![0,K_1]\!] \times ... \times [\![0,K_d]\!]} \Kg (\locpix,c_{k,i}(0)) \; \alpha_{k,i}(0),
    \label{eq:velocity_field_grid}
\end{align}
where $\alpha_{k,i}(0)$ is a momentum vector attached to $c_{k,i}(0)$.

Instead of optimizing $v_{0,i}$ by optimizing its momentum vectors $\alpha_{k,i}(0)$, we will optimize them in the wavelet domain under the constraints that the finer-scale coefficients are equal to $0$.

More precisely, we define
\begin{align*}
    ((a_{S,k,i})_{k},(d^{o}_{s,k,i})_{1 \leq s \leq S, o, k}) = FWT((\alpha_{k,i})_{k})
\end{align*}
where the Wavelet Transform has been applied to the momentum vectors component by component. Note that by construction,
\begin{align*}
    (\alpha_{k,i})_{k} = IWT  ((a_{S,k,i})_{k},\; (d^{o}_{s,k,i})_{1 \leq s \leq S, o, k}).
\end{align*}

Performing optimization with the wavelet coefficients $(a_{S,k,i})_{k}$ and $(d^{o}_{s,k,i})_{1 \leq s \leq S, o, k}$ instead of the momentum vectors $(\alpha_{k,i})_{k}$
amounts to switch from the description of \Cref{eq:velocity_field_grid} to
\rev{\begin{align*}
    v_{0,i}(\locpix) &= \sum \limits_{k}^{} a_{S,k,i} \; \widetilde{\phi}_{S,k}(\locpix) \\
    &+ \sum \limits_{s=1}^{S} \sum \limits_{o,k} d^{o}_{s,k,i} \; \widetilde{\psi}^{o}_{s,k}(\locpix)
\end{align*}}
where $\widetilde{\phi_{s,k}}$ and $\widetilde{{\psi}^{o}_{s,k}}$ are functions defined by replacing $\phi_{0,k'}$ by $\Kg (\cdot,c_{k'}(0))$
in Equations~\ref{eq:phidec} and~\ref{eq:psidec}:
\rev{\begin{align*}
    \widetilde{\phi_{s,k}} &= \sum_{k'} \gamma_{s,k,k'} \; \Kg (\cdot,c_{k'}(0))\\
    \widetilde{\psi^{o}_{s,k}} & = \sum_{k'} \gamma^{o}_{s,k,k'} \; \Kg (\cdot,c_{k'}(0)).
\end{align*}}

Even if $v_{0,i}$ is now defined through the vectorial wavelet coefficients $(a_{S,k,i})_{k}$ and $(d^{o}_{s,k,i})_{s, o, k}$ instead of the momentum vectors $(\alpha_{k,i})_k$, it remains a linear combination of the $K_{g}(x,c_k(0))$
so that we are still in the setting of Durrleman et al. \citep{durrleman2012}
and can rely on
\begin{align*}
    v_{t,i}(\locpix) =  \sum\limits_{k} \Kg (\locpix,c_{k,i}(t)) \; \alpha_{k,i}(t).
\end{align*}

Note that we do not use the Haar parameterization outside the initialization. Indeed, the initial grid is deformed under the action of the diffeomorphism $\flowphi_{t,i}$ when the time evolves, so that our Haar coefficients would be hard to interpret for $t>0$.

\begin{figure*}[tb]
\centering
\subfloat[Source and target images]
{\includegraphics[width=0.28\textwidth]{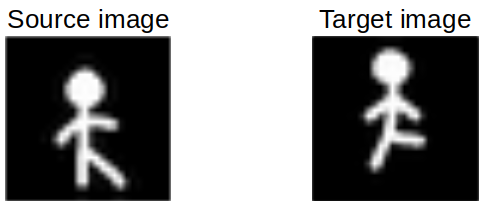}}\\
\subfloat[Classical LDDMM algorithm]{\includegraphics[width=0.43\textwidth, trim=0 -139 0 0, clip]{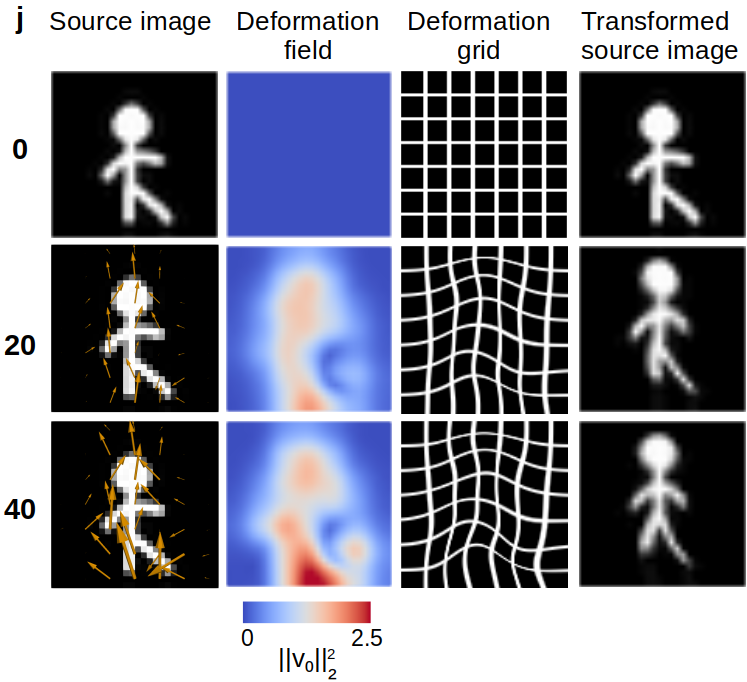}}
 \hspace{0.5cm}
\subfloat[LDDMM $+$ Multiscale algorithm]{\includegraphics[width=0.47\textwidth]{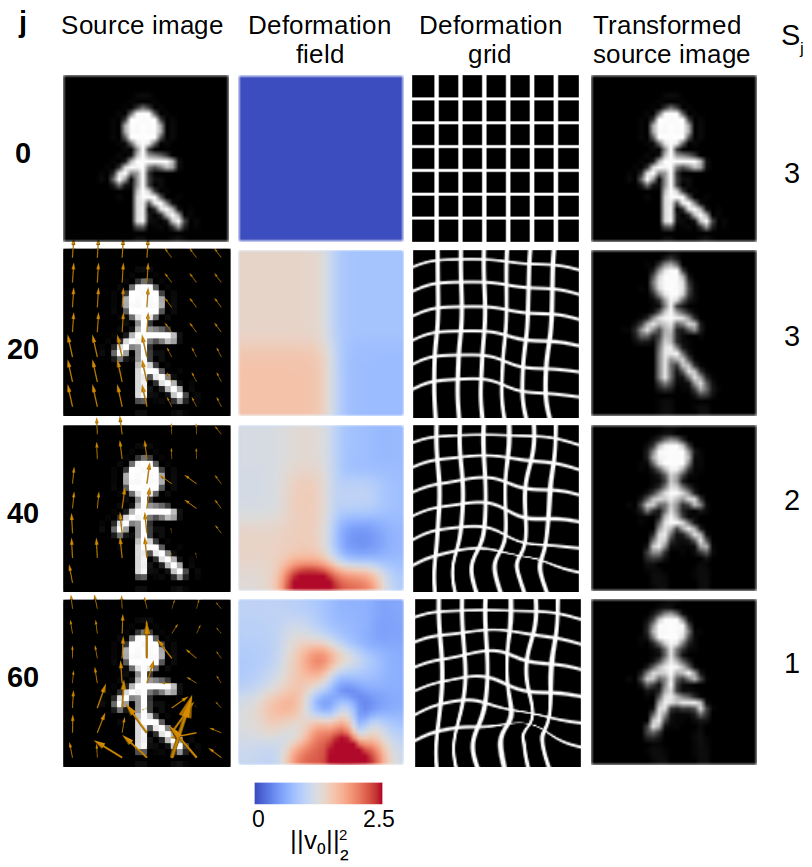}}
\caption{Original LDDMM and LDDMM combined with multiscale optimization applied to a registration example. The source and target images (panel (a)) present feature differences at a large scale (character translation) and at a finer scale (orientation of the arms and legs). Registration was performed with $\sigma_g=4$ ($k_g=49$ control points). For each algorithm, we display the source-to-target vector fields (orange arrows, scaling factor $=5$), the corresponding deformation field $L_2$ norm in RGB, the deformation grid and the transformed source image every 20 iterations until convergence. $S_j$ is the scale constraining the momentum vectors at iteration $j$. Notice how the multiscale strategy first estimates coarse displacements to move the character to the top right, then finer transformations to adjust the position of the arms. The original LDDMM algorithm applies fine-scale deformations to the entire image (except the borders), while the multiscale algorithm estimates fine deformations around the character's right leg, and smoother deformations elsewhere.}
\label{fig:example_ctf_characters}
\end{figure*}

\subsection{Coarse-to-fine atlas estimation}\label{ref:coarse_to_fine_atlas_estim} 

\subsubsection{Reparameterization of the initial velocity fields}

The key difference between our scheme and the one
of Durrleman et al. \citep{durrleman2012} is the use of the Haar parameterization in a fixed grid for the initial velocity fields. \rev{Thus, we want to optimize the following cost function:}
\begin{align}
	\label{eq:objective_function_haar}
&E(I_{ref}, (\beta_{0,i})_{ 1\leq i  \leq N}) = \nonumber\\
& \qquad\sum_{i=1}^{N} \left( \frac{d(I_i, I_{ref}  \circ \flowphi_{1,i}^{-1})^2}{2\sigma^2}  + \|v_{0,i}\|_V^2 \right),
\end{align}

where $\beta_{0,i}$ is a set of wavelet coefficients $(a_{S,k,i})_k$ and $(d^o_{s,k,i})_{s,k,o}$, related to the momentum vectors $\alpha_{0,i}$ by  $\alpha_{0,i} = IWT(\beta_{0,i}) = M_{IWT} \beta_{0,i}$. 

\rev{Since  $\alpha_{0,i} = M_{IWT} \beta_{0,i}$}, there is a relationship between the gradient of the cost function
with respect to the wavelet coefficients $\nabla_{\beta_{0,i}} E$ and the gradient with respect to the momenta $\nabla_{\alpha_{0,i}}E$:
\begin{align*}
\begin{split}
 \nabla_{\beta_{0,i}} E &= M_{IWT}^{{T}}  \nabla_{\alpha_{0,i}}E \\
    &= M_{FWT} \nabla_{\alpha_{0,i}}E \\
    &= FWT(\nabla_{\alpha_{0,i}}E)
\end{split}
\end{align*}
where we have used the fact that $M_{IWT}^{{T}}$=$M_{FWT}$ because the transform is orthonormal.

\rev{Therefore, if we compute the gradient $\nabla_{\alpha_{0,i}}E$ with Durrleman et al. numerical scheme \citep{durrleman2012}, we can then easily obtain the values of $\nabla_{\beta_{0,i}} E$ for almost the same cost as the one of this latter gradient. As illustrated in Fig. 2, this means that optimization can now be performed in the wavelet domain: we can compute $\nabla_{\alpha_{0,i}}E$ using the original LDDMM algorithm, and, instead of updating the ${\alpha_{0,i}}_i$, compute $\nabla_{\beta_{0,i}} E$ and update the $(\beta_{0,i})_i$. Without further changes, such algorithm would lead to exactly the same results as the original algorithm.}

To obtain different and hopefully better results, we enforce some constraints on the wavelet coefficients of the initial velocity fields. Namely, we use a coarse-to-fine initialization strategy by optimizing first the initial velocity fields whose wavelets coefficients are null at the finest scales and adding progressively these fine scale coefficients. 
In the following, we describe in detail our procedure which is summarized in \cref{alg:algorithm_ctf} and illustrated on a simple registration example in \Cref{fig:example_ctf_characters}. We chose a source and a target image presenting both large and small scale differences to illustrate clearly how the original and coarse-to-fine algorithms cope with multiscale deformations. 

\begin{algorithm}
\begin{algorithmic}[1]
\small
\caption{Multiscale optimization.}
\label{alg:algorithm_ctf}
\STATE{\textbf{Input}}
\STATE{Set of images ${(I_i)_{1 \leq i \leq N}}$, template image ${I_{ref}}$, kernel width  ${\sigma_g}$, 
step size $h$, \rev{initial scale $S_0$}}
\STATE{\textbf{Initialization}}
\STATE{${j \leftarrow 0}$}
\STATE{$c_0 \leftarrow$ Grid of control points with spacing ${\sigma_g}$}
\STATE{Template image ${I_{ref}(j) \leftarrow I_{ref}}$}
\STATE{Momentum vectors ${\alpha_{0,i}(j) \leftarrow 0 \quad \forall i \in [1,N]}$} 
\STATE{Initialize the gradients $\nabla_{\alpha_{0,i}} E\leftarrow 0 \quad \forall i \in [1,N]$ and $\nabla_{I_{ref}} E\leftarrow 0$}
\STATE{${\beta_{0,i}(j) \leftarrow FWT(\alpha_{0,i}(j))} \quad \forall i \in [1,N]$}
\STATE{\rev{Current scale ${S_j \leftarrow} S_0$}}
\REPEAT
\STATE{${j \leftarrow j+1}$}

\FOR{each subject $i$} 
\begin{sloppypar}\STATE{Compute the evolution of $\alpha_i(t)$ using \Cref{eq:HamSyst} \COMMENT{Gradients computation} \label{line:gradients_1}}\end{sloppypar}
\STATE{Compute $\flowphi_i$ by solving the flow equation (\Cref{eq:flow_equation})}
\STATE{Deform the template $I_{ref}$ with $\flowphi_i$}
\STATE{Compute the energy $E$ (\Cref{eq:objectiveparam})}
\STATE{Compute the gradients $\nabla_{\alpha_{0,i}} E$ {\footnotesize \&} $\nabla_{I_{ref}}E$}\label{line:gradients_2}
\STATE{Wavelet transform ${\nabla_{\alpha_{0,i}} E}$ with \cref{alg:FWT}: ${\nabla_{\beta_{0,i}} E} \leftarrow FWT(\nabla_{\alpha_{0,i}} E) = (a^i_{\smax,k})_k \cup (d^{i,o}_{s,k})_{1 \leq s \leq \smax,k,o}$}\label{line:fwt}
\FOR{each detail coefficient $d^{i,o}_{s,k}$}
\IF{${s < S_j}$ } 
    \STATE{$d^{i,o}_{s,k} \leftarrow 0$ \label{line:finer_scale_silencing_}} \COMMENT{Finer scale silencing}
\ENDIF
\ENDFOR
\ENDFOR

\STATE{${\beta_{0,i}(j) \leftarrow \beta_{0,i}(j-1) - h  \times  \nabla_{\beta_{0,i}} E}$ for all subjects i \COMMENT{Parameter update}\label{line:param_update}} 
\STATE{${\alpha_{0,i}(j) \leftarrow IWT(\beta_{0,i}(j))} \forall i \in [1,N]$ for all subjects $i$} \COMMENT{\cref{alg:IWT}}
\STATE{${I_{ref}(j) \leftarrow I_{ref}(j-1) - h  \times  \nabla_{I_{ref}} E}$}

\STATE{Compute the mean residuals $\nabla_{j}$ according to  \Cref{eq:residuals}}
\IF{${\frac{\nabla_{j-1} - \nabla_{j}}{\nabla_{j-1}} < 0.01}$ and $S_j > 1$} 
\STATE{${S_j \leftarrow S_{j-1} - 1}$ \COMMENT{Scale refinement step} \label{line:scale_refinement}}
\ENDIF

\UNTIL{Convergence}
\RETURN{Template image ${I_{ref}}$ and momentum vectors $\alpha_{0,i}$}
\end{algorithmic}
 \end{algorithm}

\subsubsection{Coarse-to-fine initialization}

The coarse-to-fine optimization
can be seen as an initialization of each new scale with the optimal template-to-subject deformations of the previous coarser scale. 
More precisely, at iteration $j$, we only optimize the wavelet coefficients of the vector fields whose scales are above or equal to a current decreasing scale $S_j$. 
Since the original RKHS setting is preserved, this can easily be done by computing the gradient $\nabla_{\alpha_{0,i}}E$ with the efficient numerical scheme of Durrleman et al. \citep{durrleman2012} (lines \ref{line:gradients_1}-\ref{line:gradients_2} in \cref{alg:algorithm_ctf}), applying FWT to $\nabla_{\alpha_{0,i}}E$ to derive the gradient with respect to $\beta_{0,i}$ (line \ref{line:fwt}) and then setting to $0$ the wavelet coefficients whose scale is strictly smaller than $S_j$ (line \ref{line:finer_scale_silencing_}).
We then update the coefficients $\beta_{0,i}$ with the modified gradient and recover the updated $\alpha_{0,i}$ using the IWT function (line \ref{line:param_update}).
If we iterate without modifying the current scale, we optimize the cost function in a subspace of functions that are simpler than in the original algorithm, the wavelet transform scale limitation acting as a regularizer.

As we want to optimize on the full set of functions defined by the momentum vectors, 
we progressively decrease the current scale $S_j$. We propose to decrease the scale when we are close to convergence at the current scale (line \ref{line:scale_refinement}). This is measured by computing the mean residual value over subjects at iteration $j$:
\begin{align}
	\Delta_{j}(x) = \frac{1}{N} \sum \limits_{i} \|I_{j,0} \circ \flowphi_{j,i}^{-1} - I_i\|_2^2,
 \label{eq:residuals}
\end{align} 
where $j$ denotes the current iteration.

If the residual decrease with respect to the previous iteration is below a threshold of $1\%$, we decide that the algorithm is close to convergence. 
In the case of \Cref{fig:example_ctf_characters}, the algorithm starts at $S_0 = 3$, performs optimization until (almost) convergence at this scale, goes to scale $2$ and performs the subsequent scale transitions in the same manner.

Our optimization procedure ensures that the momenta belonging to the same area are updated with identical values. 
 At a given scale $S_j$, the velocity fields can vary spatially only at scales coarser than $S_j-1$. In other words, at scale $S_j$, the $i^{th}$ initial velocity field implicitly writes as follows:
 \begin{align*}
v_{0,i}(\locpix)=\sum \limits_{k} a_{\smax,k,i}  \widetilde{\phi_{\smax,k}}(\locpix) \nonumber \\ 
+ \sum \limits_{s=S_j}^{\smax} \sum \limits_{k,o}{} d^o_{s,k,i} \widetilde{\psi^o_{s,k}}(\locpix) 
\end{align*}
where $\widetilde{\phi_{\smax,k}}$ and $\widetilde{\psi^o_{s,k}}$ are linear combinations of localized Gaussian kernels $\Kg (\locpix,\gk(0))$.

 


 When the algorithm reaches scale 1, the momenta are updated independently of each other. 
 Importantly, unlike previous approaches that represented deformations in a wavelet basis \citep{Amit1994, Downie1996, Wu2000, Sun2014, Gefen2004}, when the algorithm reaches this finest scale, the momentum vectors are free of constraints and the parameterization of the velocity fields is equivalent to its original definition, i.e. a sum of localized small Gaussian kernels. Thus, in theory, our coarse-to-fine algorithm could reach the same solutions as the original one, but as we will see in \cref{sec:experiments}, the coarse-to-fine numerical scheme converges to better solutions.

 \rev{The initial scale $S_0$ takes values between 1 and $\smax$, i.e. the maximum scale  of the wavelet coefficients $\beta_0,i(0)$. If $S_0=1$, the classical single-scale optimization is performed. By default, $S_0 = \smax$ so that coarse-to-fine optimization is performed from scale $\smax$ to scale 1.} 
 
 Our code is available in a public Git repository\footnote[2]{\url{https://github.com/fleurgaudfernau/Deformetrica_multiscale/}}. 

\subsubsection{Complexity}

Contrary to the previous coarse-to-fine algorithms developed in the LDDMM framework, our strategy does not add any complexity to the mathematical model. The parameterization of the velocity fields remains identical to that of Durrleman et al. \citep{durrleman2012}. Computation of the gradients and subsequent cost also remains identical. The only additional complexity arises from the algorithms FWT and IWT, which are of linear complexity.


\section{Experiments}
\label{sec:experiments}

In this section, we evaluate the benefits of the wavelet reparameterization of the initial velocity fields \rev{on different tasks and datasets. In \cref{sec:toy_experiment}, we use a registration experiment on toy data to illustrate the way our coarse-to-fine algorithm operates and we compare it to the multi-kernel algorithm \citep{Risser2011}.} 
\rev{Then, we compare the performance of our algorithm to that of the original, single-scale LDDMM version \citep{durrleman2012} on three atlas estimation tasks of increasing complexity.} The training phase consists in atlas estimation and the test phase consists in registering the estimated template image to a set of new images.

\rev{Experiments are run on an Ubuntu 18.04.5 machine equipped with a NVIDIA GPU driver with 12 GB memory. The original version of the algorithm is available in Deformetrica Version 4.3.0 \citep{Deformetricasoftware}}. Optimization for the classical LDDMM and LDDMM-multiscale algorithms relies on a gradient descent algorithm in which the step sizes $h$ are first scaled by the squared norm of the gradients and then diminished by a backtracking algorithm to guaranty a descent. \rev{Unless stated otherwise,} the following parameters are used in all experiments: $\sigma = 0.1$ in the cost criterion; initial step size $h = 0.01$; convergence threshold $= 0.0001$. The minimum number of iterations between successive coarse-to-fine steps is set to $5$. 
The initial template image for atlas estimation is given by the mean of the intensities of the training images.\\


\subsection{Toy experiment}\label{sec:toy_experiment}

\begin{figure}[ht]
\centering
\subfloat[Source image $I_s$]{\includegraphics[width=0.39\linewidth]{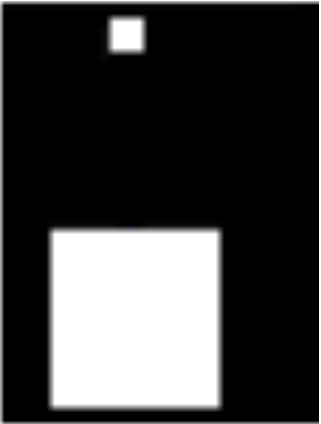}} \hspace{0.3cm}
\subfloat[Target image $I_t$]{\includegraphics[width=0.395\linewidth]{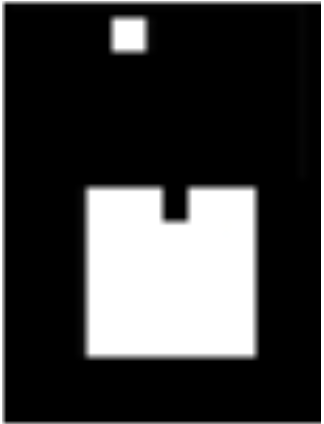}}
\caption{Registration experiment from Risser et al. \citep{Risser2011}}
\label{fig:toy_experiment}
\end{figure}

\rev{
In this experiment, we propose to assess the effect of the parameter $S_0$ on the performance of the multiscale algorithm and compare its behaviour to that of the multi-kernel framework \citep{Risser2011}. 
To this end, we reproduce a toy experiment presented in Risser et al.: the goal is to register a source image $I_s$ onto a target image $I_t$. Both images, visible in \cref{fig:toy_experiment}, contain a small (4 by 4 pixels) and a large (20 by 20 pixels) square. From $I_s$ to $I_t$, the large square is translated to the top-right and a small indentation (4 by 2 pixels) appears on the top of the square. The small square remains at the same location. Thus, the registration task requires a coarse scale deformation (displacement of the large square), a small scale deformation (creation of the indentation), all the while ensuring that the small square is not deformed. 

Performance is assessed by using the same metrics as the original experiment:
\begin{itemize}
    \item \textbf{Total residuals} $\Delta_{J}$, i.e. sum of squared differences in the source image domain at convergence
\item \textbf{Total residuals in a region of interest covering the indentation}, denoted by \DeltaJROI
\item \textbf{Standard deviation of the Jacobian determinant} (SD(J)) \citep{Leow2007}, quantifying the amount of local volume differences between the template and target images to assess the smoothness of transformations. 
\end{itemize}
In addition, we also compute the algorithms runtimes. 
}

\begin{figure*}[ht]
    \centering
    \subfloat[Residuals $\Delta_{J}$ and \DeltaJROI]{\includegraphics[width=0.28\linewidth]{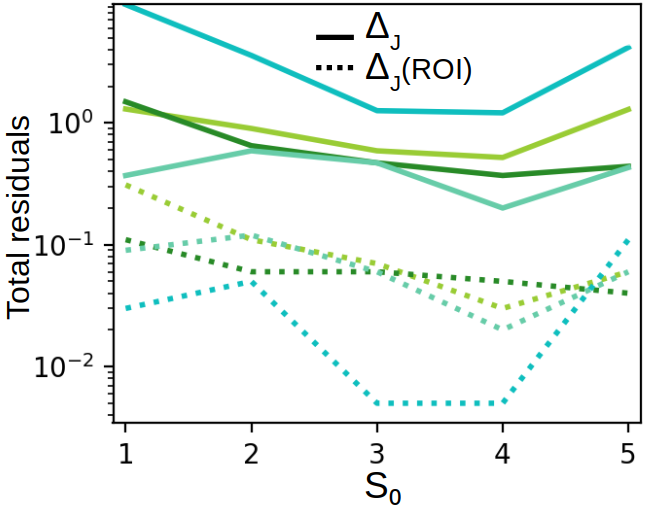}} \hspace{0.2cm}
    \subfloat[Jacobian standard deviation]{\includegraphics[width=0.283\linewidth]{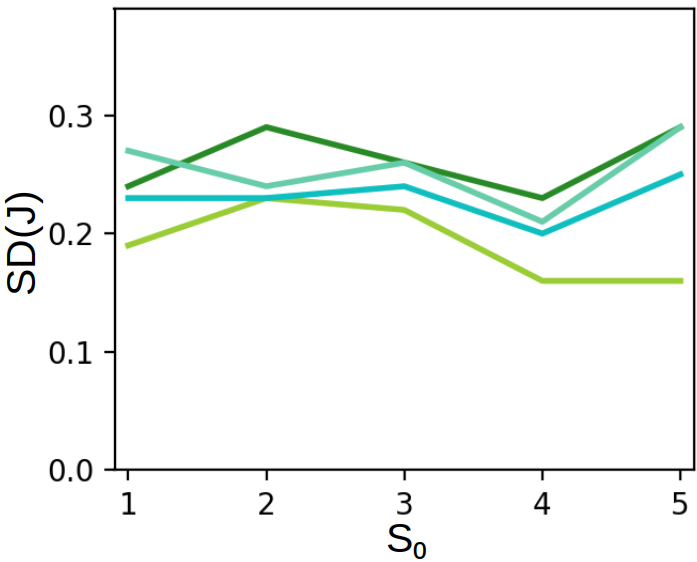}} \hspace{0.2cm}
    \subfloat[Runtime]{\includegraphics[width=0.35\linewidth]{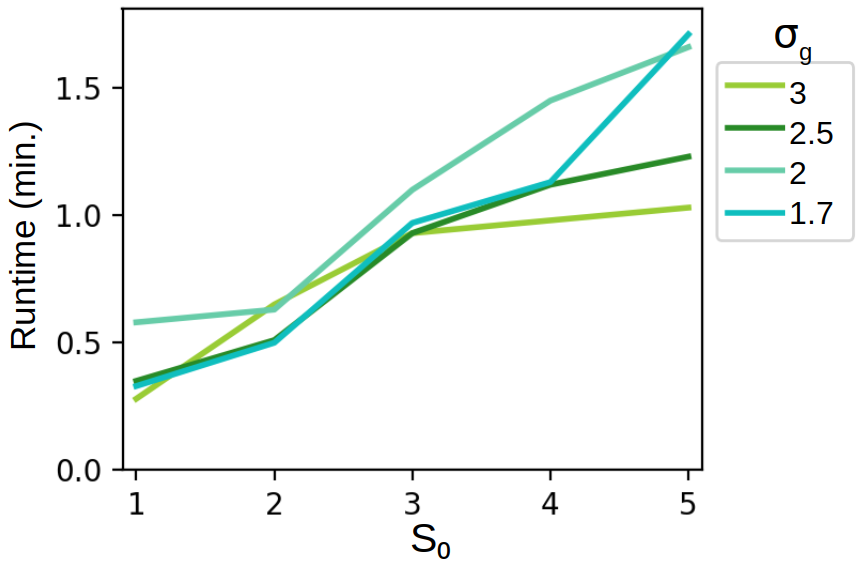}}
    \caption{Registration quality depending on the initial scale $S_0$ for different kernel widths $\sigma_g$.}
    \label{fig:scales_graphs}
\end{figure*}
\begin{figure*}[ht]
    \centering
\includegraphics[width=0.8\linewidth]{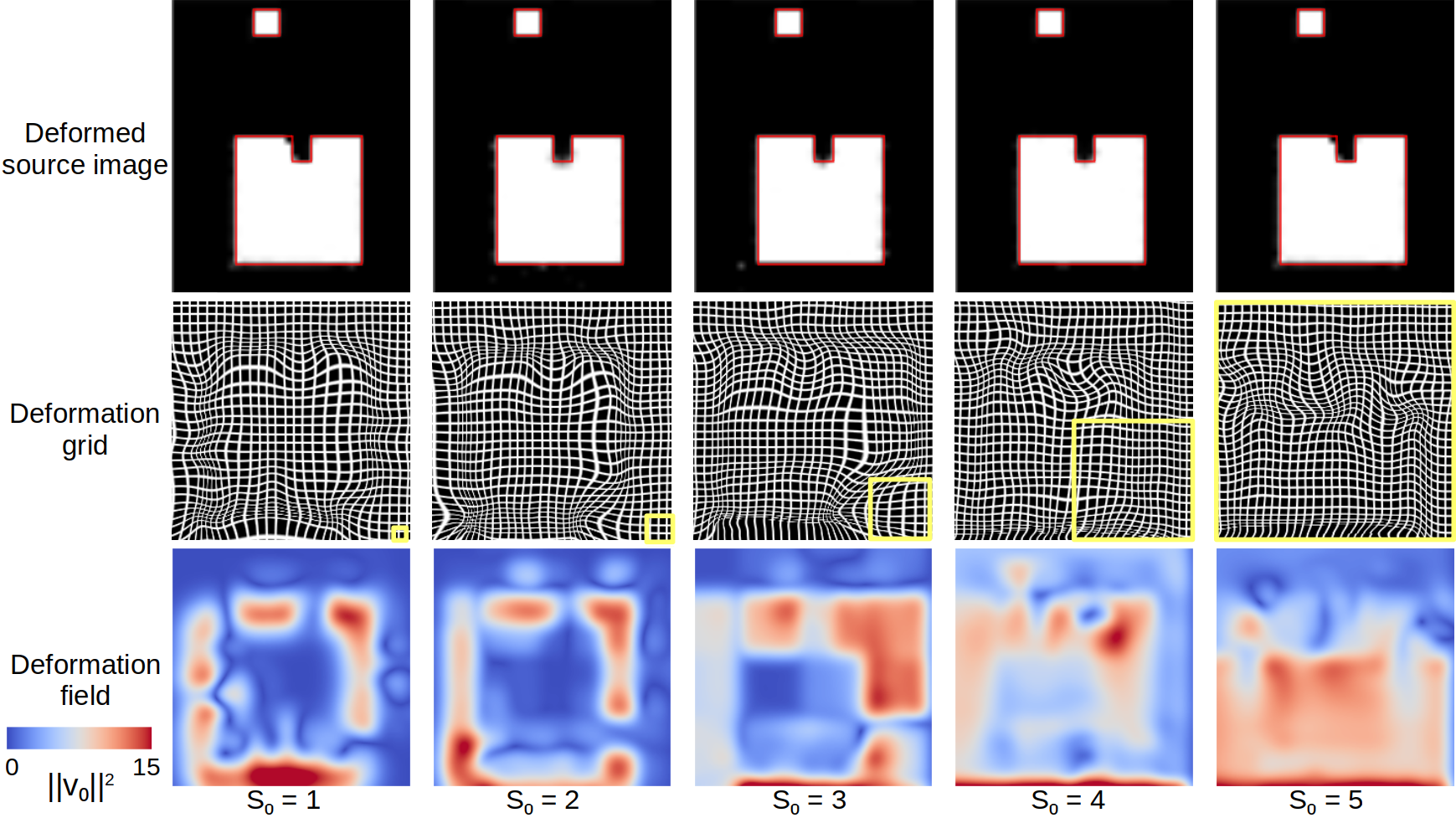}
\caption{\rev{Registration between $I_s$ and $I_t$ with the LDDMM algorithm combined with multiscale optimization using $\sigma_g=2$ and different initial scales $S_0$. Top row: deformed source images at $t=1$. The red lines represent the borders of the shapes in $I_t$. Middle row: corresponding deformation grids, zoomed in on the large square area. The yellow squares represent the support of the wavelet functions of the initial scale $S_0$, i.e. the size of the areas in which the vector fields are constant at the beginning of optimization: $32\times 32$ for $S_0=5$; $16 \times 16$ for $S_0=4$, $8\times 8$ for $S_0=3$, $4\times 4$ for $S_0=2$, and $2\times 2$ pixels for $S_0=1$. Bottom row: corresponding deformation fields $L^2$ norm in RGB, zoomed in on the large square area.}}
\label{fig:scales_grids_results}
\end{figure*}


\subsubsection{Effect of the initial scale}\label{sec:initial_scale_effect}
\rev{
In the first part of this experiment, we assess the effect of the initial scale $S_0$ on the performance of the multiscale algorithm. 
Registration between $I_s$ and $I_t$ is performed with the multiscale algorithm using kernels sizes $\sigma_g \in \{1.7, \;2, \; 2.5,\; 3\}$, corresponding to a number of control points $k_g \in \{841, \; 625,\; 400,\; 289,\}$. For each parameter $\sigma_g$, we run the multiscale algorithms with all the possible initial scales $S_0 \in \{1, 2,..., \smax\}$. Note that $S_0=1$ corresponds to the original LDDMM algorithm.

\cref{fig:scales_graphs} presents the total residuals $\Delta_{J}$ and \DeltaJROI, the standard deviation of the Jacobian SD(J) and the runtimes after registration as functions of the initial scale $S_0$. \cref{fig:scales_grids_results} displays the registration results obtained with $\sigma_g = 2$. 
Compared to our baseline (i.e. $S_0=1$), the use of our coarse-to-fine strategy decreases the total residuals value after registration. The highest performance improvement is reached for the smallest kernel $\sigma_g=1.7$: displacing the large square using fine kernels is very  costly, a problem which the multiscale strategy alleviates by imposing first large transformations. The indentation area also benefits from the multiscale strategy, as indicated by values of \DeltaJROI  inferior to 1 pixel when $S_0 > 2$. 
For all parameters $\sigma_g$, the minimum of $\Delta_J$ is reached when optimization starts at the second coarsest scale $S_0=4$. The lower performance of $S_0=5$ indicates that very coarse deformations may not benefit the optimization. At $S_0=5$, the algorithm starts by optimizing only wavelet coefficients linked to functions whose support is $[0, \; 2^4]^d$ control points. With a spacing of 2 pixels between control points, the estimated vector fields are constant on areas of $32\times 32$ pixels, which is very above the size of the large square. This effect is illustrated in \cref{fig:scales_grids_results}: with $S_0=5$, the deformation grid shows a transformation affecting a larger part of the image domain than necessary. We also notice that the final deformation fields estimated with $S_0=3$ are strongly influenced by the initial scale, again demonstration the strong impact of $S_0$ on the rest of the initialization. 

Interestingly, the multiscale strategy also seems to improve the smoothness of the final transformation when the initial scale is $S_0=4$ (\cref{fig:scales_graphs}, panel (b)). In \cref{fig:scales_grids_results}, the deformation grid obtained with $S_0=1$ displays some irregularities, whereas the deformation grid and field of $S_0=4$ are distinctively smoother. 

Finally, we observe in \cref{fig:scales_graphs}, panel (c) that the algorithm runtime increases linearly with $S_0$: for $S_0=5$, it is 3 to 4 times higher than for the original LDDMM. 

This experiment shows the importance of performing our coarse-to-fine procedure starting from a reasonably coarse scale. It can effectively enhance the outcome of the classical LDDMM when both large and fine scale transformations are needed. 

}

\subsubsection{Comparison to the multi-kernel algorithm}

\begin{table}[h]
    \caption{\rev{Performance of different algorithms on the toy experiment: original LDDMM (Deformetrica version), LDDMM combined with multiscale optimization (with $S_0$ =3 for $\sigma_g=7$ and $S_0=4$ for $\sigma_g< 7$), original LDDMM (Risser et al. version and results \citep{Risser2011}) and multi-kernel LDDMM (Risser et al. results). 
    The multi-kernel algorithm was tested with different combinations  MK\textsubscript{$n$}, with  MK\textsubscript{$n$} a sum of $n$ weighted Gaussian kernels with standard deviations linearly sampled between 1 and 10. The reported values are the total residuals values between $I_s$ and $I_t$ in the image domain and a ROI covering the indentation and standard deviation of the Jacobian determinant.}}
    \centering
\small
\begin{simpletabular}{@{}cc|ccc|ccc@{}}
\multicolumn{2}{c}{} &\multicolumn{3}{c}{\makecell{\textbf{LDDMM}\\ \textbf{(Deformetrica)}}} & \multicolumn{3}{c}{\tableLDDMMmultiscale}\\
\toprule
 $\sigma_g$ & $k_g$ & \closer{$\Delta_{J}$} & \textls[-150]{\DeltaJROI} & \close{SD(J)} & \closer{$\Delta_{J}$} & \textls[-150]{\DeltaJROI} & \close{SD(J)} \\
\midrule
7 & 64 & 11.52 & 6.67 & 0.09 & 11.06 & 6.45 &0.10 \\
3 & 289 & 1.30 & 0.31 & 0.19 & 0.52 & 0.03 & 0.16\\
2.5 & 400 & 1.49 & 0.11 & 0.24& 0.37 & 0.05 &0.23  \\
2 & 625 & 0.37 & 0.09 & 0.27 &  {0.20} & 0.02& 0.21\\
1.7 & 841 & 9.36 & 0.03 & 0.23  & 1.21 &0.0& 0.24\\ 
\bottomrule
\end{simpletabular}

\vspace{0.2cm}

\begin{simpletabular}{cccccc}
\begin{minipage}{.4\linewidth}
\begin{simpletabular}{@{}c|ccc@{}}
\multicolumn{4}{c}{\makecell{\textbf{LDDMM} \\ \textbf{(Risser et al.)}}}\\
\toprule
 $\sigma_g$ & \closer{$\Delta_{J}$} & \closer{\DeltaJROI} & \close{SD(J)} \\
\midrule
7& 8.3 & 7.2 & 0.01  \\
3& 3.0 & 2.5 & 0.02  \\
1& 0.39 & 0.07 & 0.07  \\
\bottomrule
\end{simpletabular}
\end{minipage} & & & & & 
\begin{minipage}{.7\linewidth}
\begin{simpletabular}{@{}c|ccc@{}}
\multicolumn{4}{c}{\close{\textbf{Multi-kernel LDDMM}}}\\
\toprule
 $\sigma_g$ & \closer{$\Delta_{J}$} & \closer{\DeltaJROI} & \close{SD(J)} \\
   \midrule
\close{ MK\textsubscript{3}} & 0.88 & 0.34 &0.02 \\
\close{ MK\textsubscript{4}} & 1.0 & 0.35& 0.02\\
\close{ MK\textsubscript{5}} & 0.81 & 0.34& 0.02\\
\close{ MK\textsubscript{6}} & 0.94 & 0.37& 0.02\\
\close{ MK\textsubscript{7}} & 0.94 & 0.37& 0.02\\
\bottomrule
\end{simpletabular}
\end{minipage} 
\end{simpletabular}

\label{tab:risser}
\end{table}

\rev{
In the second part of this experiment, we compare the behaviour and performance of our multiscale optimization vs the multi-kernel algorithm \citep{Risser2011}. 
Direct comparison between the two multiscale algorithms is not possible as 1) our implementation of the classical LDDMM differs from that of Risser et al. 2) experiments in Risser et al. were performed with a fixed (unknown) number of iterations, which is not possible in our setting since our multiscale strategy has an effect on the total number of iterations. 
It is however possible to compare how the multi-kernel LDDMM improves the classical LDDMM versus how our multiscale optimization improves the classical LDDMM.
To this end, we perform registration between $I_s$ and $I_t$ using the same trade-off between regularity and data attachment (i.e. 0.5)  and the same kernel widths as in the original experiment, i.e. $\sigma_g \in \{2, \; 3,\; 7\}$, plus the intermediary kernels $\sigma_g \in \{1.7, \; 2.5\}$. Contrary to Risser et al. we do not present results obtained with $\sigma_g > 7$ as the related kernels are too large for our multiscale algorithm to be of use; nor do we present results obtained with $\sigma_g=1$,  which has a low performance on this experiment with our LDDMM implementation. 
 In \cref{tab:risser}, we report the performance of the algorithms for the different values of $\sigma_g$: classical LDDMM (Deformetrica version \citep{durrleman2012}), LDDMM combined with our multiscale algorithm (with $S_0=\smax-1$), classical LDDMM (Risser et al. version and results \citep{Risser2011}), and multi-kernel LDDMM (Risser et al. results). As detailed in \cref{sec:state_of_the_art_lddmm}, flows in the multi-kernel framework are defined by a weighted sum of Gaussian kernels of different sizes: in \cref{tab:risser}, we report performance values obtained using MK\textsubscript{$n$}, $n \in \{3, 4, 5, 6, 7\}$ \citep{Risser2011}, where MK\textsubscript{$n$} denotes a sum of weighted Gaussian kernels with standard deviations linearly sampled between 1 and 10.

The classical LDDMM (Risser et al. version) reaches its best performance (in terms of residuals decrease) for the smallest kernel of size $\sigma_g=1$. Compared to this baseline, the multi-kernel LDDMM has slightly lower performance whichever combination MK\textsubscript{$n$} is used. If we compare our classical LDDMM (Deformetrica version) with our multiscale strategy, the latter yields significantly lower values of $\Delta_J$ and \DeltaJROI for all values of $\sigma_g$. Three cases are discernible: 
\begin{itemize}
    \item When a large kernel is used, i.e. $\sigma_g=7$, the original LDDMM has poor performance and the improvement yielded by our coarse-to-fine strategy is only marginal as it is constrained, at core, by the scale $\sigma_g$.
    \item As observed in \cref{sec:initial_scale_effect}, when a fine kernel is used, i.e. $\sigma_g=1.7$, the original LDDMM has poor performance. Its outcome is dramatically enhanced by the multiscale algorithm, indicating that it can successfully avoid unrealistic local minima.
    \item When an adequate kernel is used, e.g. $\sigma_g=2$, the original LDDMM has good performance (with residual values similar to that of Risser et al. baseline), and the multiscale algorithm is still capable of improving its outcome marginally, reaching the best performance out of all the compared algorithms. 
\end{itemize}

In terms of deformation smoothness, the multi-kernel LDDMM algorithm attains values of SD(J) similar to that of Risser et al. original LDDMM for $\sigma_g=3$, while reaching better performance. Compared to their baseline ($\sigma_g=1$), the multi-kernel LDDMM decreases SD(J) by a factor of 3. In contrast, our multiscale algorithm only marginally improves smoothness.


This experiment shows that our multiscale strategy can enhance the outcome of the classical LDDMM in a variety of situations, especially when using fine kernels, and can be competitive with another multiscale algorithm that builds true multiscale flows. 
In the following experiments, we will compare the performance of the classical, single-scale LDDMM algorithm with that of our multiscale algorithm on three different atlas estimation tasks. We will evaluate their ability to estimate high quality, stable template images as well as natural template-to-subjects deformations in different settings. 
The multiscale algorithm will be initialized with $S_0=\smax-1$.

}


\subsection{Handwritten digits}\label{sec:digits}

\begin{table}[h]
\caption{Performance (mean and standard deviation) of the \rev{classical LDDMM and the LDDMM combined with multiscale optimization} over five folds of cross-validation after atlas estimation and registration on the dataset of handwritten digits. The relative residual errors and SSIM are computed between each train/test image and the deformed template. Bold style indicates statistically significant differences between the algorithms ($p<0.05$).}
    \centering
\small
\vspace{-0.2cm}
\begin{tabular}{@{}cc|cc|cc@{}}
 \multicolumn{2}{c}{ $\;$} & \multicolumn{2}{c}{\textbf{R}} & \multicolumn{2}{c}{\textbf{SSIM}}\\
\toprule
 $\sigma_g$ & $k_g$ & \closer{LDDMM} & \tableLDDMMmultiscalee & \closer{LDDMM} & \tableLDDMMmultiscalee\\
 \midrule
  & &  \multicolumn{4}{c}{Atlas estimation} \\
$3$ & $100$ &
$10.2 \; \scs{(3.1)}$ & $7.9 \; \scs{(1.9)}$ & $0.74\; \scs{(0.10)}$ & $0.79\; \scs{(0.05)}$ \\
$2$ & $196$ 
& $10.1 \; \scs{(1.7)}$ &  \boldmath{$5.3\; \scs{(0.5)}$} & $0.74\; \scs{(0.05)}$ &  \boldmath{$0.79\; \scs{(0.01)}$}\\
$1.5$ & $361$ 
& $7.3\; \scs{(1.7)}$ & \boldmath{$4.3\; \scs{(1.4)}$} & $0.75\; \scs{(0.05)}$ & \boldmath{$0.79\; \scs{(0.05)}$}\\
\midrule
& & \multicolumn{4}{c}{Registration} \\
$3$ & $100$ &
$10.0 \; \scs{(6.1)}$ & $7.5 \; \scs{(6.0)}$ & $0.74\; \scs{(0.04)}$ & $0.79\; \scs{(0.03)}$ \\
$2$ & $196$    
& $8.5  \; \scs{(5.2)}$ & \boldmath{$6.5 \; \scs{(6.3)}$}& $0.72\; \scs{(0.03)}$ & \boldmath{$0.78\; \scs{(0.03)}$}  \\
$1.5$ & $361$ &     
$ 9.0\; \scs{(5.1)}$ & \boldmath{$5.3 \; \scs{(7.0)}$} & $ 0.71\; \scs{(0.04)}$ & \boldmath{$0.78\; \scs{(0.03)}$}\\
\bottomrule
\end{tabular}
\label{tab:digits_residuals}
\end{table}

\begin{figure*}[ht]
\centering
\subfloat[$\sigma_g=3; k_g=100$ (8 pixels per point)]{\includegraphics[width=0.9\textwidth]{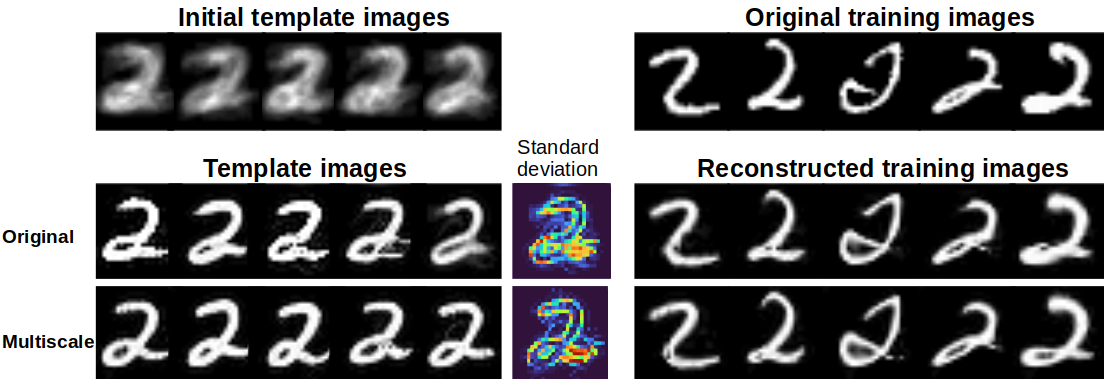}}\\
\vspace{-0.3cm}
\subfloat[$\sigma_g=2; k_g=196$ (4 pixels per point)]{\includegraphics[width=0.9\textwidth]{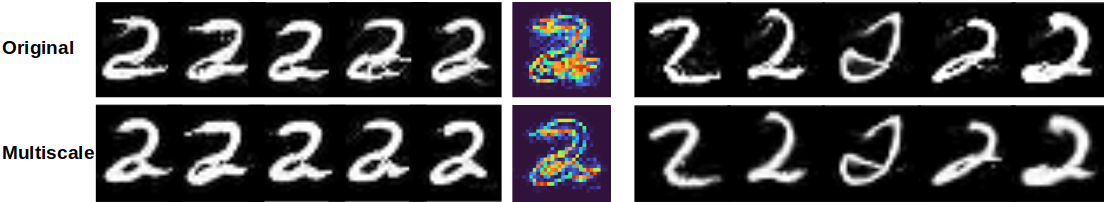}}\\
\vspace{-0.3cm}
\subfloat[$\sigma_g=1.5; k_g=361$ (2 pixels per point)]{\includegraphics[width=0.9\textwidth]{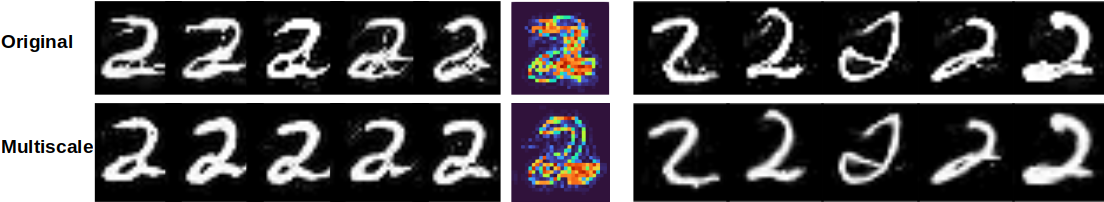}}
\caption{Estimation of the template image by the \rev{original LDDMM and the LDDMM combined with multiscale optimization} on the dataset of handwritten digits with different values of $\sigma_g$. For each experiment, five template images (estimated with non-intersecting training sets) and \rev{their standard deviation image} are presented on the left, along with the template image from the first training set warped to the first five training images on the right. $\sigma_g$: width of the Gaussian kernel; $k_g$: corresponding number of control points.} 
\label{fig:test_digits}
\end{figure*}

\rev{In this section as well as the following experiments, we evaluate the ability of our multiscale scheme to enhance the outcome of the classical LDDMM on an atlas estimation task.} We test the original \rev{LDDMM} algorithm against the coarse-to-fine strategy using an experimental procedure similar to that of Durrleman et al. \citep{durrleman2012}: 
atlas estimation is performed using 20 randomly-chosen training images and the estimated template image is registered to 10 randomly chosen test images with the same parameters as those used during training. The experiment is repeated five times with different training and test sets, with no intersection between any of the training and test sets. This procedure is  reproduced with different kernel widths $\sigma_g$. Since the five experiments are performed on independent datasets, we use paired Student t-tests to compare performance between the two algorithms.

\rev{The performance of the algorithms is assessed with the following metrics:}
\begin{itemize}
    \item \textbf{Relative residual error}: $ R = \frac{\Delta_{J}}{\Delta_0}$\\
    with $\Delta_{0}$ and $\Delta_{J}$ and  the mean residual value over subjects at iteration $0$ and at convergence, respectively.
    \item \textbf{Structural Similarity Index Metric} \citep{SSIM} (SSIM) after training and test:\\
$SSIM(I_1, I_2) = l(I_1, I_2) \times c(I_1, I_2) \times s(I_1, I_2)$\\
where $I_1$ and $I_2$ are the compared images, $l$ is a function comparing the luminance (i.e. the mean pixel/voxel intensity) of the images, $c$ compares the contrast (i.e. the standard deviation of the image intensities), and $s$ quantifies the structural similarity between $I_1$ and $I_2$ (i.e. the correlation between luminance- and contrast- normalized intensities). SSIM values range between $-1$ (dissimilarity) and $1$ (near-perfect similarity).
\item \rev{\textbf{Standard deviation of the Jacobian determinant} (SD(J)) of the displacement fields linked to the deformations $(\phi_i)_i$ to assess the regularity of the estimated transformations.} 
\item \rev{\textbf{Runtimes}}
\end{itemize}

In this section, we use images of the digit 2 extracted from the well-known United States postal database of handwritten digits \citep{Hastie2009}. The size of the images is 28 by 28 pixels. \cref{tab:digits_residuals} shows the mean relative residual error and similarity yielded by the algorithms over the experiments. \rev{Standard deviation of Jacobian values and runtime can be seen in Appendix \cref{tab:digits_jacobian}.}
\Cref{fig:test_digits} presents the five template images estimated by each algorithm with three different sets of parameters, along with the template image estimated from the first training set warped to five of the training images.

\cref{tab:digits_residuals} shows that the coarse-to-fine algorithm reaches lower residual error and higher similarity than the original \rev{LDDMM} algorithm during both atlas estimation and registration, with differences that reach significance for $\sigma_g = 2$ and $\sigma_g = 1.5$. Consistent with these results, we observe in \Cref{fig:test_digits} that the original algorithm yields highly irregular template images, a trend which worsens when the number of control points increases, indicating overfitting.
\rev{The first template image warped towards the training images yields images that are close to the original ones for $\sigma_g = 3$. However, for lower values of $\sigma_g$, the reconstructed training images are less realistic.}
These observations belie the quantitative evaluation, which shows that the performance of the original algorithm increases with the number of parameters. This discrepancy demonstrates that residual error alone is not sufficient to evaluate the accuracy of the algorithms.

Unlike the original \rev{LDDMM} algorithm, the coarse-to-fine procedure produces realistic template images that are stable across folds and whose quality is preserved when $\sigma_g$ is decreased.  
Moreover, all reconstructed images are very close to the original ones. Their quality slightly increases with the number of control points: this is most evident for the third and fifth reconstructed subjects, which become more accurate for lower values of $\sigma_g$.

\rev{Expectedly, in Appendix \cref{tab:digits_jacobian}, standard deviations of Jacobian values increase with the number of control points for both algorithms. The multiscale algorithm reaches slightly higher SD(J) values than the original algorithm, indicating more stretching of the deformed source images. 
In average, the multiscale algorithm has higher runtime, reflecting the time spent performing optimization at each scale (as shown in \cref{sec:toy_experiment}).}


\subsection{Artificial characters} \label{sec:artificials_characters}

In the previous experiment, one can remark that the performance of the algorithms diverge most when a high number of parameters is used. Therefore, one might simply be tempted to employ the original LDDMM algorithm with a lower number of parameters, as in Durrleman et al. \citep{durrleman2012} whose experiments were performed with 36 control points. 
However, datasets that present a higher amount of details and inter-subject variability may benefit from our coarse-to-fine strategy even when a lower number of control points is used. To confront our algorithm with a more difficult task, we manually designed a dataset of 30 characters. The size of the images is 28 by 28 pixels.

\begin{figure*}[ht]
\centering
\subfloat[$\sigma_g=5, k_g=36$ (22 pixels per point)]{\includegraphics[width=0.9\textwidth]{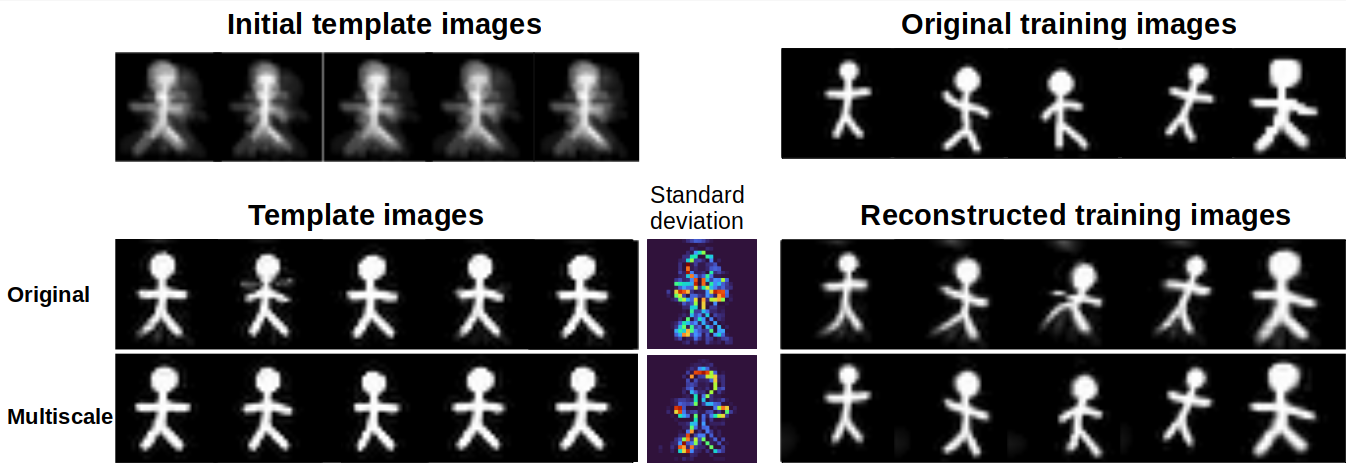}}\\ 
\vspace{-0.3cm}
\subfloat[$\sigma_g=4, k_g=49$ (16 pixels per point)]{\includegraphics[width=0.9\textwidth]{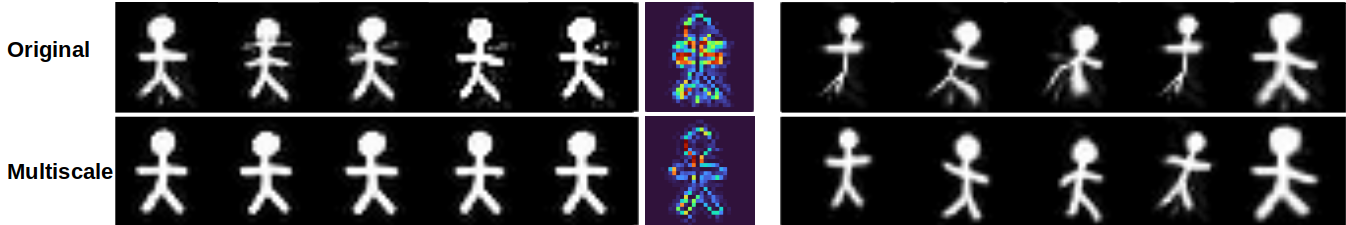}}\\
\vspace{-0.3cm}
\subfloat[$\sigma_g=3, k_g=100$ (8 pixels per point)]{\includegraphics[width=0.9\textwidth]{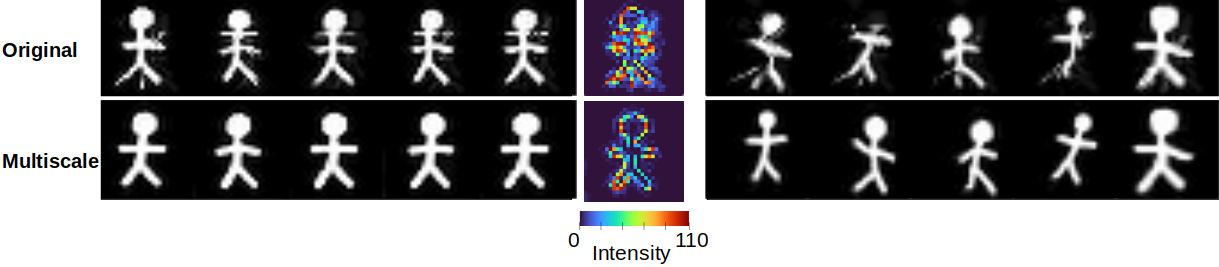}}
\caption{Estimation of the template image by the \rev{original LDDMM and the LDDMM combined with multiscale optimization} on the toy dataset with different parameters $\sigma_g$. 
For each experiment, five estimated template images (for each fold of cross-validation) \rev{and their standard deviation image} are presented on the left, along with the template image estimated from the first training set warped to the first five training images on the right.}
\label{fig:test_characters}
\end{figure*}

\begin{table}[b]
\caption{Performance (mean and standard deviation) of the \rev{classical LDDMM and the LDDMM combined with multiscale optimization} over five folds of cross-validation after atlas estimation and registration on the dataset of artificial characters. The relative residual errors and SSIM are computed between each train/test image and the deformed template.}
\centering
       \small
\begin{tabular}{@{}cc|cc|cc@{}}
\multicolumn{2}{c}{\textbf{}}& \multicolumn{2}{c}{\textbf{R}}&\multicolumn{2}{c}{\textbf{SSIM}}\\
\toprule
 $\sigma_g$ & $k_g$ & \closer{LDDMM} & \tableLDDMMmultiscalee & \closer{LDDMM} & \tableLDDMMmultiscalee\\
 \midrule
 & &  \multicolumn{4}{c}{Atlas estimation}\\ $5$ & $36$ & 
 $14.0\; \scs{(4.5)}$ & $5.5\; \scs{(0.8)}$ & $0.76\; \scs{(0.04)}$& $0.82\; \scs{(0.01)}$\\
$4$ & $49$ 
& $16.4\; \scs{(3.8)}$ & $4.4\; \scs{(0.4)}$ & $0.74\; \scs{(0.03)}$& $0.84\; \scs{(0.00)}$ \\
$3$ & $100$ 
& $18.4\; \scs{(4.5)}$ & $3.3\; \scs{(0.5)}$ & $0.68\; \scs{(0.03)}$ & $0.83\; \scs{(0.00)}$\\
\midrule
& & \multicolumn{4}{c}{Registration} \\
$5$ & $36$ & 
$9.8\; \scs{(6.2)}$ & $6.0\; \scs{(6.3)}$ & $0.75\;\; \scs{(0.02)}$& $0.78\;\; \scs{(0.03)}$ \\
$4$ & $49$ 
& $10.7\; \scs{(6.0)}$ & $4.6\; \scs{(4.2)}$ & $0.73\;\; \scs{(0.02)}$&  $0.82\;\; \scs{(0.01)}$\\
$3$ & $100$ 
& $13.4 \; \scs{(7.5)}$ & $3.6\; \scs{(3.8)}$ &$0.64\;\; \scs{(0.05)}$& $0.80\;\; \scs{(0.00)}$\\
\bottomrule
\end{tabular}
\label{tab:characters_residuals}
\end{table}

We compare our algorithm to the \rev{original LDDMM} using cross-validation: the dataset is randomly split into a training set (24 images) and a test set (6 images). Each algorithm  independently estimates  a template image from the training set, and then registers the template to each image in the test set with the same parameters as those used during training. This procedure is repeated five times, and reproduced with different values of $\sigma_g$. No statistical tests are performed because of the overlap between the training sets and between the test sets. \rev{Performance is assessed using the same metrics as in \cref{sec:digits}.} \cref{tab:characters_residuals} displays the mean relative residual error and SSIM after atlas estimation and registration and \Cref{fig:test_characters} shows the five template images estimated by each algorithm with three different sets of parameters, along with five reconstructed training images. 

\cref{tab:characters_residuals} shows that the coarse-to-fine algorithm reaches lower residual error than the original LDDMM algorithm. The performance of the coarse-to-fine strategy increases with the number of control points during training and test, while the original LDDMM algorithm demonstrates the opposite trend. These results are supported by the qualitative evaluation of the template images. 
In \Cref{fig:test_characters}, for $\sigma_g=5$, the two algorithms generate template images that present discrete but noticeable differences. With the original LDDMM version, the arms and legs of the characters appear slightly fuzzier, and the second template image is noisy. The reconstructed images yielded by the original algorithm are blurry (and even erroneous in case of the third subject), while the template and reconstructed images yielded by the coarse-to-fine algorithm seem sharp and accurate.

As in the previous experiment, the quality of the template images estimated by the original LDDMM algorithm deteriorates when the number of control points increases: with $\sigma_g = 4$ and $\sigma_g = 3$, images become fuzzier and display erroneous features inherited from the initial templates (mean intensity images), indicating high dependency on the initialization.
The morphology of all but one reconstructed characters is also completely erroneous. In contrast, the multiscale algorithm is able to produce stable, sharp and correct template images for all parameters. Similar observations can be made regarding the transformation of the template image towards the five training images: the multiscale strategy succeeds in generating images that are nearly identical to the original ones. 

These differences have a simple explanation: the original version simultaneously estimates the overall shape of the characters and details such as the location and orientation of the arms and legs, making it more dependent on the initial template image and leading to the selection of erroneous features, while the coarse-to-fine strategy first focuses on estimating the characters main features, which are then refined when the finer scales are optimized. 
This phenomenon is illustrated by movies showing template optimization across iterations, available at the first author's webpage\footnote[3]{\url{https://fleurgaudfernau.github.io/Multiscale_atlas_estimation}}.

\rev{In Appendix \cref{tab:characters_jacobian}, it can be seen that the multiscale algorithm has slightly lower SD(J) values during registration with $\sigma_g=3$ compared to the original LDDMM algorithm. This suggests a better ability of our multiscale strategy to regularize deformations with a high number of parameters.}

This experiment shows that on a dataset with high variability, the original LDDMM algorithm is unable to estimate templates that are satisfying with respect to quantitative and qualitative criteria. The coarse-to-fine algorithm outperforms the former in both criteria, regardless of the number of parameters.

\subsection{Fetal brain images} \label{sec:brain_mris}

\begin{figure*}[ht]
\subfloat[$\sigma_g=7$; $k_g=4,050$ (311 voxels per point)]{{\includegraphics[width=\linewidth]{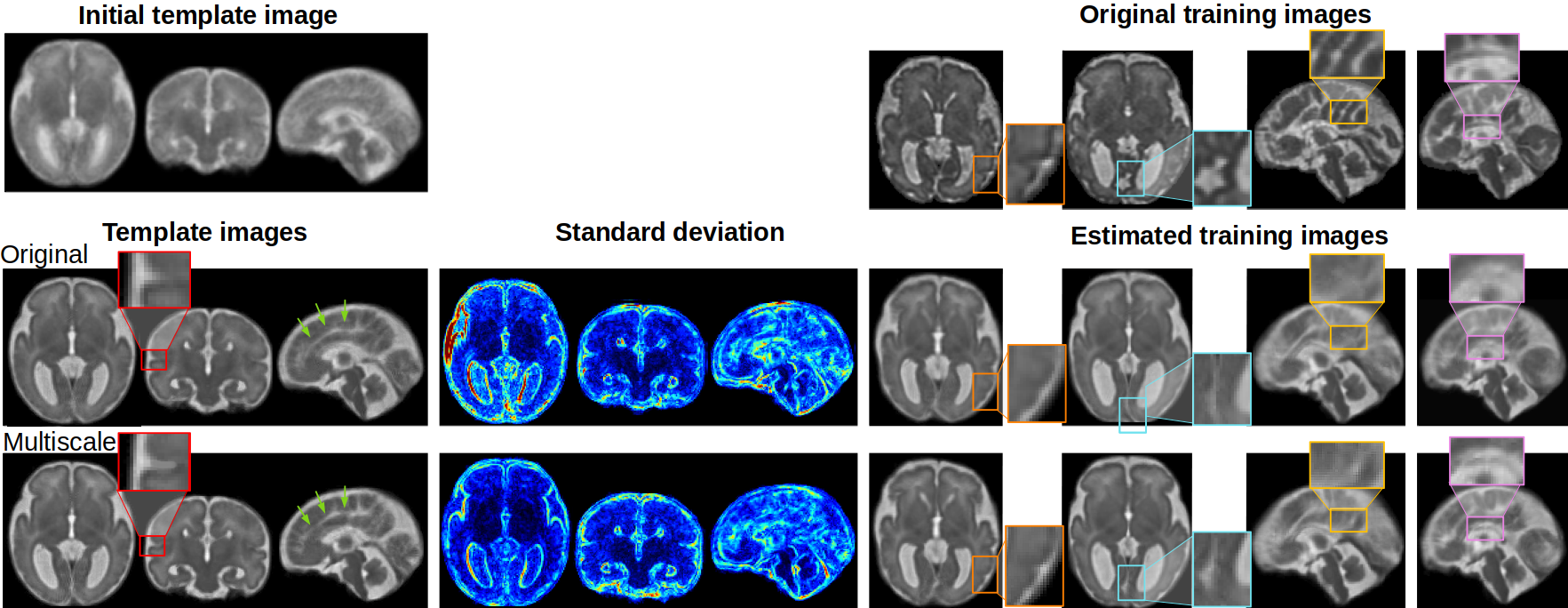}}}
\vspace{-0.41cm}
\newline
\vspace{-0.39cm}
\subfloat[$\sigma_g=5$; $k_g=10,080$ (125 voxels per point)]{
{\includegraphics[width=\linewidth]{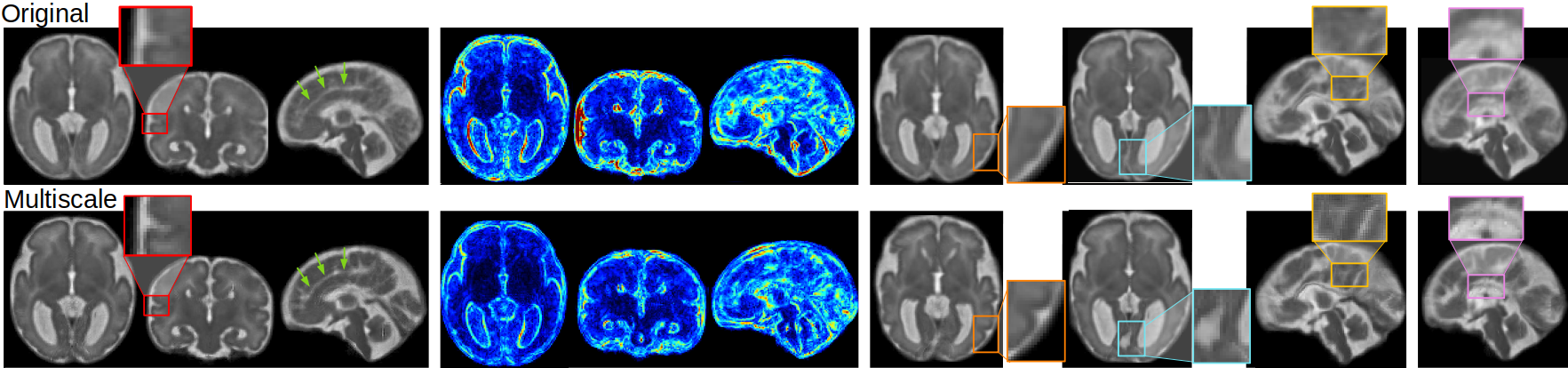}}}
\newline
\subfloat[$\sigma_g=4$; $k_g=20,250$ (63 voxels per point)]{
{\includegraphics[width=\linewidth]{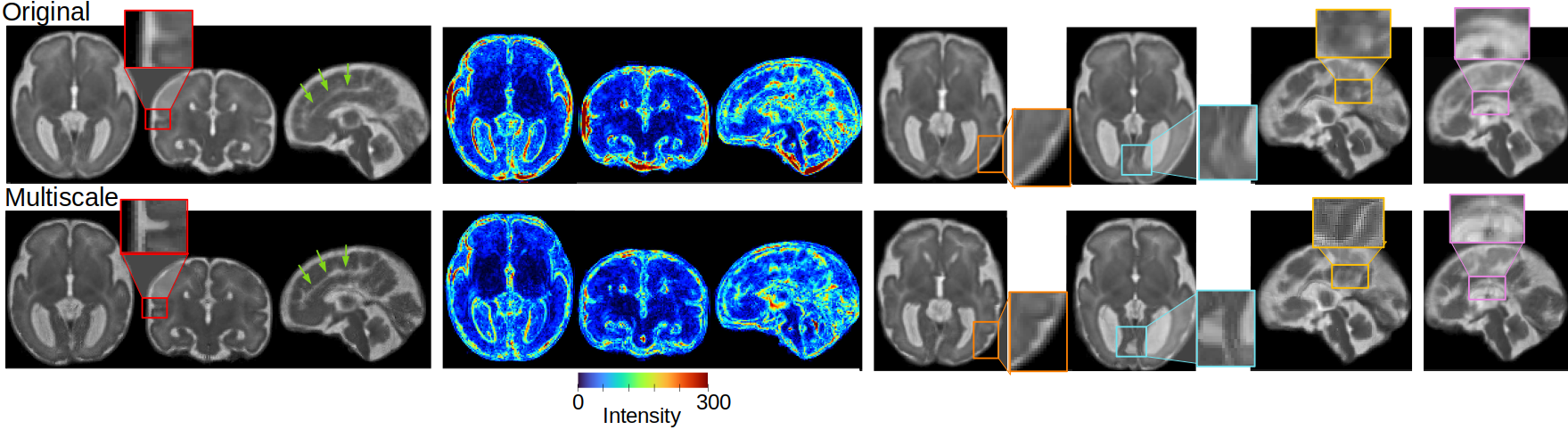}}}
\caption{Atlas estimation by the \rev{original LDDMM and the LDDMM combined with multiscale optimization} on the dataset of fetal brain images.
For each experiment, the estimated template volume from the first fold of cross-validation is presented in the left column. \rev{The first column displays the standard deviation of the template intensities over five folds of cross-validation.} The first row displays four different training images in axial or sagittal view; the four rightmost columns display the corresponding reconstructed images (i.e. the template image warped to the training images). 
Salient differences between images are zoomed in with a factor 2. Red square: superior temporal sulcus. Green arrows: cingulate sulcus. Orange square: temporal cortex. Blue square: interhemispheric fissure. Yellow square: interhemispheric sulci. Pink square: pericallosal area.}
\label{fig:test_brains}
\end{figure*}

\begin{table}[h]
    \caption{Performance (mean and standard deviation) of the \rev{classical LDDMM and the LDDMM combined with multiscale optimization} over five folds of cross-validation after atlas estimation and registration on the dataset of fetal brains images. The relative residual error (R) and the SSIM are computed between each train/test image and the deformed template.}
    \centering
\small
\hskip -1cm
\begin{tabular}{@{}cc|cc|cc@{}}
\multicolumn{2}{c}{\textbf{}} & \multicolumn{2}{c}{\textbf{R}} & \multicolumn{2}{c}{\textbf{SSIM}}\\
\toprule
 $\sigma_g$ & $k_g$ & \closer{LDDMM} & \tableLDDMMmultiscalee & \closer{LDDMM} & \tableLDDMMmultiscalee\\
\midrule
 & &  \multicolumn{4}{c}{Atlas estimation}\\ 
$7$ &  $4,050$ & $43.6 \; \scs{(3.7)}$ & $37.8\;\; \scs{(0.9)}$   & $0.73\; \scs{(0.01)}$ & $0.74 \; \scs{(0.02)}$ \\
$5$ & $10,080$  & $38.5 \; \scs{(4.0)}$ & $32.5\; \scs{(1.6)}$ & $0.74 \; \scs{(0.02)}$ & $0.77  \; \scs{(0.0)}$\\
$4$ & $20,250$ & $34.5\; \scs{(3.8)}$ & $31.9\; \scs{(3.8)}$ & $0.76\; \scs{(0.02)}$ & $0.78\; \scs{(0.02)}$ \\
 \midrule
& & \multicolumn{4}{c}{Registration} \\
$7$ &  $4,050$ & $35.0\; \scs{(6.2)}$ & $31.1\; \scs{(5.7)}$ & $0.71\; \scs{(0.01)}$ & $0.72 \; \scs{(0.0)}$\\
$5$ & $10,080$ & $31.4\; \scs{(5.0)}$ & $25.6\; \scs{(5.2)}$& $0.74 \; \scs{(0.01)}$& $0.76 \; \scs{(0.01)}$ \\
$4$ & $20,250$ & $26.7\; \scs{(1.4)}$ & $22.8\; \scs{(1.1)}$ & $0.75\; \scs{(0.01)}$ & $0.78\; \scs{(0.01)}$\\
\bottomrule
\end{tabular}
\label{tab:brains_residuals}
\end{table}

To evaluate the performance of our coarse-to-fine approach on a dataset of clinical images, we use 30 fetal brain MRIs with agenesis of the corpus callosum acquired in Hopital Trousseau, France \citep{gaudfernau2021}. Gestational ages are comprised between 32 and 34 weeks of gestation (mean = 32.9 ${\pm}$ 0.6).

Agenesis of the corpus callosum is a developmental anomaly characterized by the total or partial absence of the corpus callosum. It is often associated to anatomical features such as widening of the lateral ventricles. Atlas estimation can help better understand congenital anomalies by providing an insight into how these anatomical characteristics vary together \citep{gaudfernau2021}. However, as abnormal fetal brains may present a wide range of defects, this makes atlas estimation more difficult and prone to errors than with datasets of healthy fetuses. Thus, it is crucial to develop algorithms that are able to estimate realistic templates on both healthy and abnormal subjects.

The brain MRIs are preprocessed, volume reconstructed and rigidly aligned according to the procedure described in Gaudfernau et al. \citep{gaudfernau2021}. The final images have size 105x100x120 voxels. 
Cross-validation is performed in the same manner as in \cref{sec:artificials_characters}, with 24 images used for atlas estimation and 6 images used for testing. \Cref{fig:test_brains} presents an example of estimated template image during cross validation along with 4 reconstructed training images for different parameters $\sigma_g$, and \cref{tab:brains_residuals} displays the mean relative residual error and SSIM after atlas estimation (training) and registration (test). Visual examination of the templates is performed by an expert radiologist. 

\cref{tab:brains_residuals} shows that the coarse-to-fine algorithm achieves lower residual error and higher SSIM values than the original algorithm during atlas estimation and registration for all values of $\sigma_g$. \rev{In Appendix \cref{tab:brains_jacobian}, both algorithms show similar values of SD(J), though slightly higher for the multiscale algorithm, suggesting a better ability to capture some of the volume changes characterizing abnormal subjects, e.g. shrunk corpus callosum. }
The template images in \Cref{fig:test_brains} present subtle but noticeable differences: the brain volumes estimated with the multiscale optimization display sharper features and enhanced contrast between structures, especially for higher values of $\sigma_g$. The template images estimated by the original LDDMM algorithm display more fuzzy areas, especially at the boundary between cortical gray matter and white matter (see for example the superior temporal sulcus, red squares) and the medial surface of the brain (e.g. the cingulate sulcus, green arrows). Compared to the original algorithm, the multiscale strategy yields template images that are more stable across folds, especially in regions with high inter-subjects variability, e.g. cortical folds, lateral ventricles and corpus callosum area. 

The multiscale template-to-subject transformations also build more accurate anatomical structures: see for example the more pronounced gyration patterns (orange squares) and  clearly delineated interhemispheric fissure and lateral ventricle (blue squares). Interestingly, the multiscale template image warped to the training subjects reveals more abnormal features associated to corpus callosum agenesis, such as the typical radiating sulci (yellow squares) and missing corpus callosum (pink squares).  

\rev{As in the previous experiments, the multiscale algorithm has higher runtimes than the original LDDMM algorithm in Appendix \cref{tab:brains_jacobian}.}

Altogether, these results indicate that our coarse-to-fine strategy can successfully enhance the results of atlas estimation applied to real-world, complex clinical data.

\section{Discussion} 
\label{sec:discussion}

\rev{In this paper, we took advantage of the hierarchical property of the wavelet decomposition to develop a coarse-to-fine optimization procedure in the LDDMM framework. Specifically, we proposed a Haar-like wavelet representation of the initial velocity fields to enhance the outcomes of the classical RKHS-based LDDMM algorithm. The transfer of information from coarse to fine scales ensures smarter initialization of the deformations at each level, leading the algorithm to favor more accurate solutions and avoid unrealistic local minima. }
Contrary to previous coarse-to-fine algorithms introduced in the LDDMM framework \citep{Sommer2013, Gris2018, Modin2018, Miller2020}, our approach adds no complexity to the mathematical model. 
Specifically, the reparameterization of the velocity fields can be seen as an additional layer of spatial regularization, which preserves the RKHS structure of the vector fields and the efficient numerical scheme used to compute the gradients. 
This reparameterization can easily be translated to other mathematical frameworks which model deformations using vector fields. \rev{For example, our multiscale strategy could be combined with models that express diffeomorphisms as a composition of deformations of increasingly fine scales \citep{Modin2018, Miller2020}.}

\rev{We first performed a registration experiment to assess the influence of the initial scale parameter and compare our algorithm to the multi-kernel strategy \citep{Risser2011}. Results suggest that our multiscale algorithm can be competitive with the multi-kernel framework and lead to more accurate matching. 
While the focus of the multi-kernel framework is the design of an explicitly multiscale deformation model, our strategy's primary goal is efficiency, in the spirit of classical coarse-to-fine strategies. A clear advantage of our multiscale optimization is its convenience: we only need to identify \textit{one} relatively well-performing kernel and apply our multiscale scheme, while in the multi-kernel framework, one needs to identify \textit{several} relevant kernels, combine them, and optimize their weights, with no guarantee to converge to a better local minima \citep{Risser2011}. 
However, our method does not improve deformation smoothness, which reflects the different bases upon which the two algorithms are built: our strategy is multiscale in the sense that it solves the optimization problem in spaces of increasing resolution; the multi-kernel framework is multiscale in the way deformations are defined as coexisting flows. Thus, it would be interesting to adapt our multiscale scheme to enhance optimization in the multi-kernel framework in order to estimate true multiscale flows and favor deformation smoothness.}

\rev{We then evaluated the performance of our multiscale algorithm on three atlas estimation tasks of increasing difficulty.} Compared to the original LDDMM algorithm, the coarse-to-fine algorithm yields higher quality templates with better stability \rev{(estimated by the variability of the template intensities across cross-validation folds)}, that are able to generalize to unseen images. Not only does our strategy produce images that have a realistic anatomy, but it leads to enhanced preservation of anatomical details, including unusual or abnormal ones. This makes it particularly appropriate for tasks involving high inter-subject variability, specifically clinical images. \rev{Our results suggest that the multiscale LDDMM algorithm can estimate a more diverse range of transformations all the while preserving reasonable smoothness.}

Some limitations of this algorithm have to be highlighted. 
Unlike approaches based on a mixture of kernels, our deformation model relies on a single Gaussian kernel. 
While this provides the advantage of introducing no additional parameters, the results of our algorithm depend on the choice of the kernel width $\sigma_g$ - but to a lesser extent than the original version. \rev{As shown in \cref{sec:digits}, the multiscale \rev{LDDMM} algorithm fails to avoid irregular deformations when using a very small kernel. In such scenarios, a solution could reside in local adaptation, i.e. maintain smoothness constraints in areas where the matching is near-perfect and constrain the velocity fields to be unevenly smooth, in line with strategies based on simultaneously coexisting flows \citep{Risser2011}.}
\rev{As illustrated in \cref{sec:toy_experiment}, the runtime of our multiscale optimization increases with the number of scales of the coarse-to-fine procedure. This can be prohibitive when working with a high number of parameters, which increases the number of scales, or, as shown in \cref{sec:brain_mris}, when performing atlas estimation on clinical images.}
Methodological improvements could be made to our multiscale algorithm. \rev{The Haar wavelet, which produces sharp transitions between nearby areas of the vector fields, could be replaced with smoother wavelet functions such as the Daubechies wavelet to favor more regular vector fields.}

The properties of the wavelet transform offer interesting avenues to explore. Notably, the wavelet functions $\phi_{s,k}$ and  $\psi^o_{s,k}$ are normalized depending on their scale $s$. By modifying the weights attributed to the fine and coarse scales in this normalization, one can change the relative importance attributed to high and low frequency coefficients during optimization. It would be interesting to refine our wavelet-based spatial regularizer in this way and observe how it impacts the results.

The simplicity of our algorithm makes it easy to implement, which opens up interesting perspectives. The coarse-to-fine strategy could be applied to other types of atlases such as spatio-temporal ones \citep{Debavelaere2020}, or to other statistical frameworks such as the Bayesian Mixed Effect Model \citep{Allassoniere2015}. In regard to the latter point, a significant advance would be to integrate our reparameterization of the velocity fields into this Bayesian framework by introducing priors on the wavelet coefficients, in the spirit of Downie et al. \citep{Downie1996} who decomposed deformations into a Haar wavelet basis and modelled the coefficients as independent random variables with a mixture distribution. 


In addition, we will focus on developing a dual coarse-to-fine strategy, by applying a hierarchical representation to the images, as already attempted in other mathematical models, e.g. with B-spline deformations \citep{ Rueckert1999, Loeckx2007, staring2009} and in the hyperelasticity framework \citep{debroux2021}. Alternating both coarse-to-fine strategies would very likely provide template images of even higher quality. \rev{Ultimately, we will also apply our multiscale strategies to other models available in the LDDMM framework, namely geodesic regression  \citep{Fletcher2011} and its variant piecewise geodesic regression \citep{Chevallier2017} and evaluate their usefulness on complex clinical challenges such as the modelling of the fetal brain growth in a continuous manner \citep{Licandro2016}, a well-knowingly challenging task in Computational Anatomy.}
\newline

\subsection*{Declarations}

{\small
\textbf{Acknowledgments.}
 This work was partly funded by the third author's chair in the PRAIRIE institute funded by the French national agency ANR as part of the "Investissements d'avenir" programme under the reference ANR-19- P3IA-0001.

\textbf{Conflict of interest.} The authors have no conflicts of interest to declare that are relevant to the content of this article.

\textbf{Ethics approval.} This work follows appropriate ethical standards in conducting research and writing the manuscript, following all applicable laws and regulations regarding treatment of human subjects. The experiments conducted on the dataset of fetal brain MRI has been granted ethics approval by the Institutional Review Board of the Comité d’éthique de la recherche en imagerie médicale (CERIM) under the reference CRM-2112-215.}
\newline

{\small \textbf{Data availability}
 The dataset of handwritten digits is publicly available at: \url{https://www.kaggle.com/datasets/bistaumanga/usps-dataset}. The dataset of artificial characters is available at the first author's webpage: \url{https://fleurgaudfernau.github.io/Multiscale_atlas_estimation/}. The dataset of fetal brain MRI is not available due to ethical considerations. }
\newline

{\small \textbf{Author's contributions.}
 Conceptualization: F.G.; Methodology: F. G., E.L.P.; Validation: E.B.; Writing - original draft preparation: F.G.; Writing - review and editing: E.L.P., E.B., S.A.; Funding acquisition: S.A.; Resources: E.B.; Supervision: E.L.P., S.A.}





\begin{appendices}

\begin{figure*}[ht]
\centering
\includegraphics[width=0.85\linewidth]{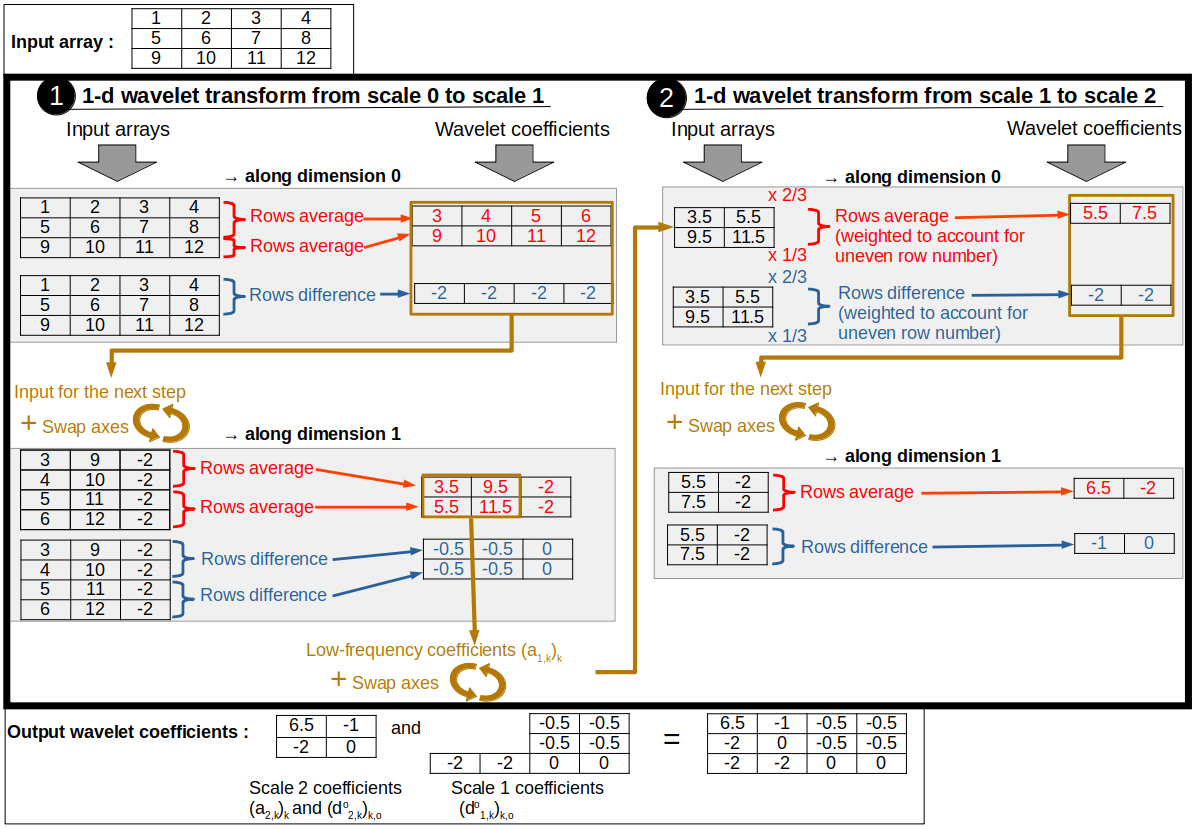}
    \caption{The FWT Algorithm (\cref{alg:FWT}) applied to a 3-by-4 array. Steps in grey squares illustrate the 1D Haar Forward step (\cref{alg:FWT_1d}). Running the IWT Algorithm (\cref{alg:IWT}) amounts to running the illustrated steps backward. For the sake of simplicity, the renormalization step (line \ref{lst:line:renorm} in \cref{alg:FWT}) is not featured. }
    \label{fig:algorithm_wavelet}
\end{figure*}

{\section{Algorithms FWT and IWT} \label{sec:algorithm_wavelets}}

To compute the Forward and Inverse Wavelet Transforms, explicit computation of the related matrices $\mathbf{M_{FWT}}$ and $\mathbf{M_{IWT}}$ is not necessary since the algorithms rely on local operations on basis coefficients. 
To compute the wavelet coefficients of an array $\boldmath{X}$ of dimension $d$, the FWT algorithm (\cref{alg:FWT}) relies on a 1D Haar Forward algorithm (\cref{alg:FWT_1d}), which computes local means and differences along one axis. Uneven numbers of rows/columns are handled by computing \textit{weighted} averages and differences so that the boundaries are given the same importance as the rest of the array. \cref{fig:algorithm_wavelet} illustrates how the FWT algorithm operates on a 2-dimensional array, corresponding to  lines \ref{line:init}-\ref{line:end1} in \cref{alg:FWT}.

Finally, the output wavelet coefficients have to be normalized to preserve the input signal energy: this is done by computing the matrix $\mathbf{M_{FWT}}$ in \cref{alg:renorm} and computing the renormalization matrix $\mathbf{R}$ from $\mathbf{M_{FWT}}$:
\begin{align*}
    \mathbf{R}[i]=\frac{1}{\|\mathbf{M_{FWT}}[i,:]\|_2}
\end{align*}
where $\|\mathbf{M_{FWT}}[i,:]\|_2$ is the $L_2$ norm of the $i_{th}$ row of $\mathbf{M_{FWT}}$.

The IWT algorithm (\cref{alg:IWT}) runs in a manner that is symmetrical to the FWT algorithm, by using the 1D Haar Backward algorithm (\cref{alg:IWT_1d}) to compute finer-scale coefficients from coarse scale coefficients one axis at a time.

\begin{algorithm}[H]
\begin{algorithmic}[1]
\small
\caption{1D Haar Forward step}
\label{alg:FWT_1d}
\STATE{\textbf{Input}}
\STATE{$\beta$: array of shape $(k_1, ...,k_D)$, $d$: axis along which to compute the transform, $K_d$: size of the original array along axis $d$,
$s$: current scale, $w_s$: list of scales for each axis $[1,...,D]$ in ascending order}
\STATE{\textbf{Initialization}}
\STATE{Swap axes in $\beta$ to put axis $d$ in position $0$}
\IF{$k_d > 1$} 
    \IF{$k_d$ is even} 
        \STATE{$\beta_{a} \leftarrow (\beta[0::2] + \beta[1::2]) / 2$} \COMMENT{Average the consecutive rows of $\beta$}
        \IF{$K_d \neq {2 \times k_d}^s$} 
            \STATE{$\delta \leftarrow K_d/2^s - (k_d - 1)$} \COMMENT{Weighting of the border}
             \STATE{$\beta_{a}[-1] \leftarrow (\beta[-2] + \delta * \beta[-1]) / ( 1+ \delta)$}
          \ENDIF
          \STATE{$\beta_{d} \leftarrow \beta[0::2] - \beta_{a}$} \COMMENT{Difference between the consecutive rows of $\beta$}
          
    \ELSE
        \STATE {$\beta_{a} \leftarrow (\beta[0:-1:2] + \beta[1::2]) / 2$} \COMMENT{Average the consecutive rows of $\beta$}
         \STATE {$\beta_{d} \leftarrow \beta[0:-1:2] - \beta_{a}$} \COMMENT{Difference between the consecutive rows of $\beta$}
          \STATE {$\beta_{a}  \leftarrow$ concatenate$([\beta_{a}, \beta[-1]])$}  \COMMENT{Add the last unpaired row}   
    \ENDIF
    \STATE{$\beta \leftarrow $ concatenate$([\beta_{a}, \beta_{d}])$}

\ENDIF
 \STATE{$w_s[d] \leftarrow [\ceil{k_d/2}] + w_s[d]$}
\RETURN{$\beta, w_s$}
\end{algorithmic}
 \end{algorithm}

 \begin{algorithm}[H]
\begin{algorithmic}[1]
\small
\caption{1D Haar Backward step}
\label{alg:IWT_1d}

\STATE{\textbf{Input}}
\STATE{$\beta$: array of wavelet coefficients of shape $(k_1, ..., k_D)$, $d$: axis of $\beta$ along which to apply the backward transform, $K_d$: size of the original array along axis $d$, $s$: current scale, $w_s$: list of scales for each axis $[1,...,D]$ in ascending order}

\STATE{\textbf{Initialization}}
\STATE{$\beta \leftarrow$ Swap axes in $\beta$ to put axis $d$ in position $0$}
\STATE{$X \leftarrow$ zero array of the shape of $\beta$}
\STATE{$n_{low} \leftarrow w_s[d][0]$}
\STATE{$\beta_a \leftarrow \beta[:n_{low}]$} \COMMENT{Low-frequency coefficients}
\STATE{$\beta_d \leftarrow \beta[n_{low}:]$} \COMMENT{High-frequency coefficients}


\IF{$w_s[d][1] >1$} 
    \IF{$w_s[d][1]$ is even} 
    \STATE{$X[0::2] \leftarrow \beta_a + \beta_d$}
    \STATE{$X[1::2] \leftarrow 2 \times \beta_a - X[0::2]$}    
        \IF{$K_d \neq {2 \times k_d}^s$} 
            \STATE{$\delta \leftarrow K_d/2^s - (k_d - 1)$}
             \STATE{$X[-1] \leftarrow ((1 + \delta) \times \beta_a[-1] - X[-2])/ \delta$}
          \ENDIF
          
    \ELSE
    \STATE{$X[0:-1:2] \leftarrow \beta_a[:-1] + \beta_d$}
    \STATE{$X[-1] \leftarrow \beta_a[-1]$}  
     \STATE{$X[1::2] \leftarrow 2 \times \beta_a[:-1] - X[0:-1:2]$}   
    \ENDIF
\ENDIF
\STATE{$w_s[d] \leftarrow w_s[d][1:]$}
\RETURN{$X, w_s$}
\end{algorithmic}
 \end{algorithm}

\begin{algorithm}[H]
\begin{algorithmic}[1]
\small
\caption{Renormalization}
\label{alg:renorm}
\STATE{ } \COMMENT{Compute the Haar Forward Matrix $M_{FWT}$ and the renormalization matrix $R$}
\STATE{\textbf{Input}}
\STATE{$(K_1, ..., K_D)$: shape of the array to wavelet transform}
\STATE{\textbf{Initialization}}

\STATE{$n \leftarrow \prod \limits_{i=1}^{D} K_i$}
\STATE{$M_{FWT} \leftarrow$ zero array of shape $(n,n)$}
\FOR{$i$ in range($n$)}
\STATE{$Z \leftarrow$ zero array of shape $(K_1, ..., K_D)$}
\STATE{$i_{th}$ element of $Z \leftarrow$ 1}
\STATE{$\beta, w_s, R \leftarrow$ \textbf{Haar\_Forward}(Z, $\rho=0$)} \COMMENT{\cref{alg:FWT}}
\STATE{$M_{FWT}[:,i] \leftarrow $ flatten($\beta$)}

\ENDFOR
\STATE{$R \leftarrow$ zero array of shape $(K_1, ..., K_D)$}
\FOR{$i$ in range($n$)}
\STATE{$i^{th}$ element of $R \leftarrow \frac{1}{\|M_{FWT}[i,:]\|_2}$}
\ENDFOR
\RETURN{$R$}
\end{algorithmic}
 \end{algorithm}

 \begin{algorithm}[H]
\begin{algorithmic}[1]
\small
\caption{Forward Haar Wavelet Transform (FWT algorithm)}
\label{alg:FWT}
\STATE{ } \COMMENT{Compute the wavelet coefficients $\beta$ of an array $X$}
\STATE{\textbf{Input}}
\STATE{$X$: array of shape $(K_1, ..., K_D)$, $\rho$: renormalization factor (default: 1)}
\STATE{\textbf{Initialization}} \label{line:init}
\STATE{$\beta \leftarrow X$} \COMMENT{Array to store the wavelet coefficients}
\STATE{$w_s \leftarrow [[K_1], ..., [K_D]]$} \COMMENT{For each axis, store the scales for which the wavelet coefficients have been computed in ascending order}
\STATE{$S_{max} \leftarrow \ceil{\log_2(\max(K_1, ..., K_D)}$}
\FOR{$s$ in range($S_{max}$)}
\STATE{$\beta_{current} \leftarrow \beta[:w_s[0][0], ..., :w_s[D][0],]$} \COMMENT{Low-frequency coefficients at scale $s-1$}

\FOR{$d$ in range($D$)}
\STATE{$\beta_{current}, w_s, R \leftarrow$ \textbf{Haar\_forward\_1d\_step}($\beta_{current}$,$d$, $K_d$, $s$, $w_s$)}
\COMMENT{\cref{alg:FWT_1d}}
\STATE{$\beta[:w_s[0][0], ..., :w_s[D][0],] \leftarrow \beta_{current}$} \label{line:end1} 
\ENDFOR
\ENDFOR
\IF{$\rho \neq 0$}
\STATE{$R \leftarrow$ \textbf{renormalize}$((K_1, ..., K_D))$} \COMMENT{\cref{alg:renorm}} \label{lst:line:renorm}
\STATE{$\beta \leftarrow \beta \times R^\rho$} 
\ENDIF
\RETURN $\beta, w_s, R$
\end{algorithmic}
 \end{algorithm}

  \begin{algorithm}[H]
\begin{algorithmic}[1]
\small
\caption{Inverse Haar Wavelet Transform (IWT algorithm)}
\label{alg:IWT}
\STATE{ } \COMMENT{Given its wavelet coefficients $\beta$, compute array $X$ (i.e. coefficients of scale $0$)}
\STATE{\textbf{Input}}
\STATE{$\beta$: array of wavelet coefficients of shape $(K_1, ..., K_D)$, $\rho$: renormalization factor (default: 1), $w_s$: list of scales for each axis $[1,...,D]$, stored in ascending order}
\STATE{\textbf{Initialization}}
\STATE{$X \leftarrow \beta$}
\STATE{$S_{max} \leftarrow \ceil{\log_2(\max(K_1, ..., K_D)}$}
\IF{$\rho \neq 0$}
\STATE{$R \leftarrow$ \textbf{renormalize}$((K_1, ..., K_D))$} \COMMENT{\cref{alg:renorm}}
\STATE{$X \leftarrow X / R^\rho$} 
\ENDIF
\FOR{$s$ in range($S_{max}-1, -1, -1$)}
\STATE{$X_{current} \leftarrow X[:w_s[0][1],..., :w_s[D][1]]$}
\FOR{$d$ in range($D$)}
\STATE{$X_{current} \leftarrow$ \textbf{Haar\_Backward\_1d\_step}($X_{current}$,$d$, $K_d$, $s$, $w_s$)} \COMMENT{\cref{alg:IWT_1d}}
\STATE{$X[:w_s[0][1],..., :w_s[D][1]] \leftarrow X_{current}$}
\ENDFOR
\ENDFOR
\RETURN{$X$}

\end{algorithmic}
 \end{algorithm}

\section{Algorithms smoothness and runtime}

\rev{
\vspace{-0.5cm}
\begin{table}[h]
    \caption{\rev{Standard deviation of the Jacobian determinant (SD(J)) and runtimes (RT) (mean and standard deviation) of the classical LDDMM and the LDDMM combined with multiscale optimization over five folds of cross-validation after atlas estimation and registration on the dataset of handwritten digits. Bold style indicate statistically significant differences between the algorithms ($p<0.05$).}}
    \centering
\small
\begin{tabular}{@{}cc|cc|cc@{}}
\multicolumn{2}{c}{\textbf{}} & \multicolumn{2}{c}{\textbf{SD(J)}} & \multicolumn{2}{c}{\textbf{RT (min.)}}\\
\toprule
 $\sigma_g$ & $k_g$ & \closer{LDDMM}  & \tableLDDMMmultiscalee & \closer{LDDMM}  & \tableLDDMMmultiscalee\\
 \midrule
 & &  \multicolumn{4}{c}{Atlas estimation}\\
 $3$ & $100$ &
$0.15 \; \scs{(0.02)}$ & $0.19 \; \scs{(0.01)}$  &$8.50 \; \scs{(3.26)}$ & $13.22 \; \scs{(2.88)}$ \\

$2$ & $196$ 
& $0.18 \; \scs{(0.01)}$ &  $0.21 \; \scs{(0.01)}$ & $8.61 \; \scs{(1.29)}$ &  \boldmath{$12.29 \; \scs{(1.13)}$} \\

$1.5$ & $361$ & $0.20 \; \scs{(0.01)}$ & $0.22 \; \scs{(0.01)}$ & $7.0 \; \scs{(1.44)}$ & \boldmath{$12.34 \; \scs{(1.12)}$}\\
\midrule
& & \multicolumn{4}{c}{Registration} \\
$3$ & $100$ &
$0.20 \; \scs{(0.01)}$ & $0.20 \; \scs{(0.01)}$ &
$0.31 \; \scs{(0.05)}$ & \boldmath{$0.52 \; \scs{(0.03)}$} \\

$2$ & $196$    
& $0.22 \; \scs{(0.01)}$ & $0.24 \; \scs{(0.02)}$ & $0.38 \; \scs{(0.05)}$ & \boldmath{$0.56 \; \scs{(0.02)}$} \\
$1.5$ & $361$ & $0.23 \; \scs{(0.02)}$ & $0.26 \; \scs{(0.03)}$ & $0.32 \; \scs{(0.03)}$ & \boldmath{$0.58 \; \scs{(0.08)}$} \\
\bottomrule
\end{tabular}
\label{tab:digits_jacobian}
\end{table}

\vspace{-0.5cm}
\begin{table}[t]
\caption{\rev{Standard deviation of the Jacobian determinant (SD(J)) and runtimes (RT) (mean and standard deviation) of the classical LDDMM and the LDDMM combined with multiscale optimization over five folds of cross-validation after atlas estimation and registration on the dataset of artificial characters.}}
\centering
\small
\begin{tabular}{@{}cc|cc|cc@{}}
\multicolumn{2}{c}{\textbf{}} & \multicolumn{2}{c}{\textbf{SD(J)}} & \multicolumn{2}{c}{\textbf{RT (min.)}}\\
\toprule
 $\sigma_g$ & $k_g$ & \closer{LDDMM} & \tableLDDMMmultiscalee & \closer{LDDMM} & \tableLDDMMmultiscalee\\
 \midrule
 & &  \multicolumn{4}{c}{Atlas estimation}\\ 
 $5$ & $36$ & 
 $0.10 \; \scs{(0.01)}$ & $0.11 \; \scs{(0.01)}$ & $9.51 \; \scs{(3.04)}$ & $13.65 \; \scs{(2.52)}$ \\
 
$4$ & $49$ 
& $0.12 \; \scs{(0.01)}$ & $0.13 \; \scs{(0.01)}$ & $9.43 \; \scs{(2.19)} $ & $13.78 \; \scs{(2.62)}$  \\

$3$ & $100$ 
& $0.13 \; \scs{(0.01)}$& $0.12 \; \scs{(0.01)}$ & $11.10 \; \scs{(3.77)}$& $13.21 \; \scs{(1.77)}$ \\
\midrule
& & \multicolumn{4}{c}{Registration} \\
$5$ & $36$ & 
$0.13 \; \scs{(0.02)}$ & $0.12 \; \scs{(0.01)}$  & $0.37 \; \scs{(0.04)}$ & $0.56 \; \scs{(0.05)}$\\

$4$ & $49$ 
& $0.14 \; \scs{(0.01)}$ &  $0.14 \; \scs{(0.01)}$ & $0.38 \; \scs{(0.03)}$ &  $0.51 \; \scs{(0.05)}$  \\

$3$ & $100$ 
& $0.17 \; \scs{(0.02)}$ & $0.14 \; \scs{(0.01)}$& $0.39 \; \scs{(0.03)}$ & $0.47 \; \scs{(0.02)}$\\
\bottomrule
\end{tabular}
\label{tab:characters_jacobian}
\end{table}

\vspace{-0.5cm}
\begin{table}[t]
    \caption{\rev{Standard deviation of the Jacobian determinant (SD(J)) and runtimes (RT) (mean and standard deviation) of the classical LDDMM and the LDDMM combined with multiscale optimization over five folds of cross-validation after atlas estimation and registration on the dataset of fetal brains images.}}
    \centering
\small
\begin{tabular}{@{}cc|cc|cc@{}}
\multicolumn{2}{c}{\textbf{}} & \multicolumn{2}{c}{\textbf{SD(J)}} &\multicolumn{2}{c}{\textbf{RT (min.)}}\\
\toprule
 $\sigma_g$ & $k_g$ & \closer{LDDMM} & \tableLDDMMmultiscalee &\closer{LDDMM} & \tableLDDMMmultiscalee\\
\midrule
 & &  \multicolumn{4}{c}{Atlas estimation}\\ 
$7$ &  $4,050$    & $0.02 \; \scs{(0.01)}$ & $0.03  \; \scs{(0.01)}$  &
$134.40 \; \scs{(20.67)}$ & $189.97 \; \scs{(67.39)}$ \\

$5$ & $10,080$ & $0.02 \; \scs{( 0.01)}$ & $0.03  \; \scs{( 0.01)}$ & $95.80 \; \scs{(34.61)}$ &  $191.37 \; \scs{(89.82)}$\\
$4$ & $20,250$ & $0.03 \; \scs{(0.01)}$ & $0.05 \; \scs{(0.01)}$ & $111.38 \; \scs{(37.24)}$ & $189.74 \; \scs{(77.0)}$ \\
 \midrule
& & \multicolumn{4}{c}{Registration} \\
$7$ &  $4,050$& $0.03 \; \scs{(0.01)}$ & $0.03  \; \scs{(0.01)}$ & $4.12 \; \scs{(0.53)}$ & $6.66 \; \scs{(0.62)}$ \\

$5$ & $10,080$ & $0.03 \; \scs{( 0.01)}$ & $0.03 \; \scs{( 0.01)}$ & $4.27 \; \scs{(0.48)}$ & $7.22 \; \scs{(0.42)}$\\

$4$ & $20,250$ & $0.04 \; \scs{(0.01)}$ & $0.05 \; \scs{(0.01)}$ & $3.68 \; \scs{(0.70)}$ & $7.84 \; \scs{(0.14)}$\\
\bottomrule
\end{tabular}
\label{tab:brains_jacobian}
\end{table}
}
\end{appendices}

\clearpage

\setlength{\bibsep}{0pt}
\footnotesize
\bibliography{sn-bibliography}
\end{document}